\documentclass [a4paper] {article}

\usepackage[utf8]{inputenc}

\usepackage[french]{babel}

\usepackage{amssymb,amsmath, amsfonts}
\usepackage{mathabx}

\usepackage [dvips] {graphicx,color}

\usepackage{amsthm}

\usepackage{MnSymbol}

\newcommand*{\twosquig}{%
\mathrel{\vcenter{\offinterlineskip
\hbox{$\rightsquigarrow$}\vskip-.08ex\hbox{$\rightsquigarrow$}}}}

\usepackage{cjhebrew}

\usepackage{accents}

\newtheorem{prop} {Proposition} 
 
\newtheorem{thm} [prop] {Théorème}

\theoremstyle{definition}
\newtheorem{df}{Définition} 
\newtheorem*{df*}{Définition}

\theoremstyle{remark}
\newtheorem{rmq}{Remarque} 
\newtheorem{example}{Exemple} 
\newtheorem{exm}[example]{Exemple}

\author{Stéphane \textsc{Dugowson}
\footnote{Laboratoire Quartz / Supmeca. Email : s.dugowson@gmail.com}
}

\title {Dynamiques en interaction \\ \ (une introduction à la théorie des\\ dynamiques sous-fonctorielles ouvertes)}
\date{29 août 2016} 

\begin{document}

\maketitle
%
%

\paragraph{Résumé.}
 La théorie des dynamiques sous-fonctorielles ouvertes est une nouvelle théorie qui permet de définir des systèmes dynamiques généralisés en interaction, ces interactions produisant de nouvelles dynamiques susceptibles, bien entendu, d'entrer à leur tour dans d'autres interactions. Une grande partie du contenu de cet article se trouve déjà dans les deux textes disponibles en ligne mais non publiés \cite{Dugowson:20150807} et \cite{Dugowson:20150809}, et cela a été partiellement exposé dans les conférences \cite{Dugowson:20160202} et \cite{Dugowson:20160525}. Toutefois, il était nécessaire de donner une nouvelle présentation, unifiée et donc plus commode à consulter, de ce matériau, et de l'illustrer de quelques exemples. Dans cet article, nous introduisons en outre 
les notions nouvelles d'interaction normale et d'interaction concrète, et nous remplaçons les synchronisations \og rigides\fg\, considérées jusque là, par des synchronisations \og souples\fg,  considérablement plus générales. \`{A} noter que ce qui depuis 2011  était appelé%
\footnote{Voir notamment \cite{Dugowson:201112},\cite{Dugowson:201203},\cite{Dugowson:20150807},\cite{Dugowson:20150809}.} %
 \og dynamiques catégoriques\fg\, (respectivement \og dynamiques sous-catégoriques\fg) sera désormais désigné comme \emph{dynamiques fonctorielles} (respectivement \emph{dynamiques sous-fonctorielles}).
 
\subparagraph{\emph{Mots clés.}} Interaction. Systèmes ouverts. Dynamiques. Structures connectives.

\paragraph{Abstract.} \textsc{Interacting Dynamics (Introduction To The Theory Of Open Sub-Functorial Dynamics)---} 
The \emph{theory of open sub-functorial dynamics} is a new theory that defines interacting generalized dynamical systems. The  interactions between these dynamics produce new dynamics which, of course, can then enter into other interactions. A major part of this article can already be found in two unpublished texts \cite {Dugowson:20150807} and \cite {Dugowson:20150809} and it has been partially exposed in conferences \cite{Dugowson:20160202} and \cite {Dugowson:20160525}. However, we need to give a new, unified and therefore more convenient presentation of this material, and we also need some examples to illustrate it. Moreover, we introduce in this article the new concepts of ``normal interaction" and ``concrete interaction", and replace the previously used rigid synchronizations by much more general flexible ones.



\section*{Introduction}

La théorie mathématique que nous présentons ici est de nature systémique, en ce qu'elle permet d'élaborer des systèmes de systèmes (etc.) à partir de trois ingrédients de base : 
\begin{itemize}
\item d'abord, ce que nous appelons les \emph{dynamiques sous-fonctorielles ouvertes}, présentées dans la section 
section \textbf{§\,\ref{sec DSFO}},
qui constituent une généralisation considérable des systèmes dynamiques classiques à la fois parce qu'elles ne sont pas déterministes en général, parce qu'elles reposent sur une temporalité qui n'est pas nécessairement linéaire et parce qu'elles sont \og ouvertes\fg\, en ce sens qu'elles peuvent entrer, \emph{via} une paramétrisation, en interaction avec d'autres dynamiques, 
\item ensuite, les \emph{interactions} entre dynamiques (section \textbf{§\,\ref{subs interactions dans une famille}}), qui sont définies comme des relations entre les réalisations (les trajectoires) et les paramètres des dynamiques en jeu,
\item enfin, les \emph{synchronisations} (section \textbf{§\,\ref{subs synchronisation}}), où la temporalité de chaque dynamique en jeu est référée à celle d'un chef d'orchestre.
\end{itemize}

Ces trois ingrédients permettent de constituer ce que nous appelons des \emph{familles interactives} (section \textbf{§\,\ref{subs familles interactives}}). Une telle famille interactive donne alors naissance (section \textbf{§\,\ref{sec engendrement dynamique}}) à de nouvelles dynamiques sous-fonctorielles ouvertes  qui sont des sortes de synthèses de toutes les dynamiques en interaction dans la famille considérée. La clé de voûte de cet engendrement est le \emph{théorème de stabilité sous-fonctorielle} (théorème \ref{thm stabilite sous fonctorielle}, section \textbf{§\,\ref{subs thm stab}}) qui affirme que les dynamiques en jeu étant sous-fonctorielles, les dynamiques produites le seront encore
\footnote{Le sens du mot \emph{stabilité} dans l'expression \og théorème de stabilité sous-foncto\-rielle\fg\, n'a donc rien à voir avec celui qu'il a dans l'expression \og stabilité des systèmes dynamiques\fg. Nous l'employons ici uniquement pour souligner la stabilité \emph{conceptuelle} de la notion même de dynamique sous-fonctorielle.
}.
Du reste, c'est précisément parce que les dynamiques \emph{fonctorielles} ouvertes, dont la définition plus simple nous avait d'abord incité à les prendre comme ingrédients de base de notre théorie de l'interaction, produisent parfois dans leurs interactions autre chose que des dynamiques fonctorielles (exemples \ref{exm une dynamique sur R non reguliere} et \ref{exm autre dynamique non reguliere}, section \textbf{§\,\ref{subs exemple engendrement}}) que nous avons dû chercher à établir un cadre plus large que celui des seules dynamiques sous-fonctorielles.

Tout au long de cet article, nous suivrons en particulier l'exemple de la famille interactive que nous appellerons $\mathbb{WHY}$ ou, en caractères hébraïques, $\textcjheb{why}$, 
définissant d'abord chacune des trois dynamiques en jeu dans cette famille (sections \textbf{§\,\ref{subsubs exm futur}}, \textbf{§\,\ref{subsubs exm passe}} et \textbf{§\,\ref{subsubs exm intemporel}}), puis la famille interactive elle-même (section \textbf{§\,\ref{subsubs famille why}}) et enfin les dynamiques engendrées par cette famille, en particulier celle que nous notons $\mathbb{S}$ 
ou par la lettre hébraïque \emph{shin}, \textcjheb{/s} ((section \textbf{§\,\ref{subsubs exemple shin}}). Cet exemple nous a été inspiré par la théorie du philosophe Pierre-Michel Klein sur ce qu'il appelle la \emph{métachronologie} \cite{KleinPM:2014}.\\

Nous commençons par préciser nos notations et faire quelques rappels. Dans la version actuelle, il y a une table des matières à la fin du texte.

\subsection{Notations et rappels}

\subsubsection{Réels et intervalles}\label{subsubs reels}

\begin{itemize}
\item  $\mathbf{R}$ désigne l'ensemble des réels,
$\mathbf{R}_+$ désigne l'ensemble des réels  positifs ou nuls,
$\mathbf{R}_+^*$ l'ensemble des réels strictement positifs, $\overline{\mathbf{R}}_+=[0,+\infty]$, etc. 
\item $\mathcal{I}_\mathbf{R}$ désigne l'ensemble des intervalles ouverts de $\mathbf{R}$,
\item pour toute partie $A\subset \mathbf{R}$, on notera $\mathring{A}$ ou $int(A)$ l'intérieur de $A$, et $\overline{A}$ son adhérence,
\item pour tout $(a,b)\in \mathbf{R}^2$, on pose $]a,b[=\{t\in\mathbf{R}, a<t<b\}$.
\end{itemize}


\subsubsection{Classes de fonctions numériques}\label{subsubs notations classes applications}

\begin{itemize}
\item Pour tout intervalle $I\subset\mathbf{R}$ (non nécessairement ouvert), nous posons \[\mathcal{C}(I)=\{f:I\rightarrow \mathbf{R}, f\mathrm{\,est\,continue\,sur\,} I\}, \]
et, plus généralement, pour tout $k\in\mathbf{N}\cup\{+\infty\}$, 
\[\mathcal{C}^k(I)=\{f:I\rightarrow \mathbf{R}, f\mathrm{\,est\,de\,classe\,}\mathcal{C}^k\mathrm{\,sur\,} I\}.\]
\end{itemize}

\begin{rmq}\label{rmq Ck(a,b) avec a>b}
Pour tout $k\in\mathbf{N}\cup\{+\infty\}$,  l'ensemble $\mathcal{C}^k(\emptyset)$ est un singleton, dont l'unique élément, qui sera également noté $\emptyset$, est l'inclusion canonique $\emptyset:\emptyset\hookrightarrow\mathbf{R}$. En particulier, 
\[b\leq a \Rightarrow \mathcal{C}^k(]a,b[)=\{\emptyset\}\neq \emptyset.\]
\end{rmq}

Nous posons en outre
\begin{itemize}
\item $\mathcal{C}=\bigcup_{I\in{\mathcal{I}_\mathbf{R}}}
\mathcal{C}(I)$,
\item $\mathcal{C}^k=\bigcup_{I\in{\mathcal{I}_\mathbf{R}}}
\mathcal{C}^k(I)$,
\item $\mathcal{C}^k_{\vartriangleright}=\bigcup_{r\in ]0,+\infty]}\mathcal{C}^k(]-\infty,r[)$.
\end{itemize}

\begin{rmq}
$\emptyset\notin \mathcal{C}^k_{\vartriangleright}$.
\end{rmq}

Pour tout intervalle $K\subset \mathbf{R}$, ouvert ou non, on note
$Lip^1(K)$ l'ensemble de toutes les fonctions définies et $1$-lipschitziennes sur $K$. On pose 
\[
Lip^1=
\bigcup_{I\in{\mathcal{I}_\mathbf{R}}} {Lip^1(I)}.
\]
Pour tout $c\in\overline{\mathbf{R}}_+=[0,+\infty]$, on note
\[Lip^1([0,c\vert)=Lip^1([0,c])\cup Lip^1([0,c[).\]
En particulier,  pour $c=+\infty$, on a $Lip^1([0,+\infty\vert)= Lip^1([0,+\infty[)$, et pour $c=0$ on a $Lip^1([0,0\vert)\simeq \mathbf{R}\cup \{\emptyset\}$.
On pose en outre
\[
Lip^1_+=
\bigcup_{c\in\overline{\mathbf{R}}_+} {Lip^1([0,c\vert)},
\]
et
\[
(Lip^1_+)^*=Lip^1_+ \setminus \{\emptyset\}.
\]

\subsubsection{Catégories et graphes}\label{subsubs cat et graph}

Dans tout l'article, $\mathbf{C}$ désigne une petite catégorie. En outre, si $\mathbf{E}$  désigne une catégorie quelconque,
\begin{itemize}
\item la classe des objets de $\mathbf{E}$ est notée $\dot{\mathbf{E}}$, 
\item la classe des flèches de  $\mathbf{E}$ est notée $\overrightarrow{\mathbf{E}}$,
\item pour tout objet $A$ de $\mathbf{E}$,  $Id_A$ désigne le morphisme identité,
\item le graphe associé à $\mathbf{E}$ en \emph{oubliant} la composition des flèches sera noté\footnote{Alors que dans \cite{Dugowson:20150807}, il était noté $Gr(\mathbf{E})$.} $\vert \mathbf{E} \vert$,
\item pour toute flèche $e:A\rightarrow B$ de $\mathbf{E}$, et plus généralement pour toute arête $e$ d'un graphe, $dom(e)$ désigne sa source (ou domaine) $A$, et $cod(e)$ désigne son but (ou codomaine) $B$,
\item un couple de flèches $(g,f)\in \overrightarrow{\mathbf{E}}^2$ est dit composable si $dom(g)=cod(f)$. 
\end{itemize}

\subsubsection{Catégorie vide}  La catégorie vide, qui n'a pas de flèche, sera noté $\emptyset$ ou $\mathbf{0}$.

\subsubsection{Monoïdes} 

Un monoïde est une catégorie ayant un unique objet, que nous noterons généralement $\bullet$. Par exemple, \og la catégorie $\mathbf{C}=\mathbf{R}_+$ \fg\, désigne la catégorie $\mathbf{C}$ définie par $\dot{\mathbf{C}}=\{\bullet\}$ et $\overrightarrow{\mathbf{C}}=(\mathbf{R}_+,+)$.

En particulier, j'appelle \emph{catégorie ponctuelle}, et je note $\mathbf{1}$ ou $\{\bullet\}$, la plus petite catégorie non vide, dont l'unique objet sera noté $\bullet$  et l'unique flèche $\mathbf{0}$.

\subsubsection{Catégorie des ensembles} 

La catégorie des ensembles est notée $\mathbf{Sets}$. Ses flèches sont les \emph{applications}. 

Pour tout ensemble $E$, $\mathcal{P}E$ désigne l'ensemble des parties de $E$, et $\mathcal{P}^*E$ l'ensemble des parties non vides de $E$. L'union disjointe de deux ensembles $U$ et $V$ est notée $U\sqcup V$. Si $U$ et $V$ sont disjoints, on prend $U\sqcup V=U\cup V$.

\subsubsection{Transitions}

Pour les transitions, nous reprenons les notations et les notions de \cite{Dugowson:201112}, \cite{Dugowson:201203}, \cite{Dugowson:20150807} et \cite{Dugowson:20150809}. En particulier :\\
\begin{itemize}
\item une transition $f:A\rightsquigarrow B$  d'un ensemble $A$ dans un ensemble $B$ est une application $f:A\rightarrow \mathcal{P}B$, où $\mathcal{P}B$ désigne l'ensemble des parties de $B$; ainsi,  pour tout $a\in A$, on a $f(a)\subset B$,
\item $\mathbf{P}$ désigne la catégorie dont les objets sont les ensembles et dont les flèches sont les  transitions, dont la composition est notée $\odot$,
\item les transitions de $A$ dans $B$ s'identifient trivialement aux relations binaires
\footnote{Sur les relations binaires, voir aussi les rappels et notations de la section \textbf{§\,\ref{subsubs relations binaires applications fonctions}}.} 
 de $A$ vers $B$, de sorte que la catégorie $\mathbf{P}$ coïncide avec la catégorie usuellement désignée comme \og catégorie des relations\fg, 
\item pour tout ensemble non vide $M$, $\mathbf{P}^{\underrightarrow{\scriptstyle{M}}}$ désigne la catégorie dont les objets sont les ensembles et dont les flèches sont les familles indexées par $M$ de transitions,
\item la composition des flèches dans $\mathbf{P}^{\underrightarrow{\scriptstyle{M}}}$ est encore notée $\odot$,
\item  si $S$ et $T$ sont deux ensembles et  $u:S\leadsto T$ et $v:S\rightsquigarrow T$ deux transitions, nous écrirons 
$u\subset v$ pour exprimer le fait que 
$\forall a\in S, u(a)\subset v(a)$, 
\end{itemize}

\begin{rmq}\label{rmq appli det}
Une transition $f:A\rightsquigarrow B$ qui vérifie $card(f(a))=1$ pour tout $a\in A$ est dite \emph{déterministe} et s'identifie trivialement à une application que nous noterons encore $f:A\rightarrow B$. Ainsi, dans ce cas, $f(a)$ désignera en général un élément de $B$, bien que selon le contexte cela puisse aussi désigner un singleton inclus dans $B$. En particulier, l'application identité $Id_A:A\rightarrow A$ définit une transition $A\rightsquigarrow A$ que nous noterons encore $Id_A$ dès lors que le contexte permettra de comprendre que, pour tout $a\in A$, l'écriture $Id_A(a)$ devra être comprise comme désignant le singleton
\[Id_A(a)=\{a\}.\]
\end{rmq}

\begin{rmq}\label{rmq fonctions quasi det} Le mot \emph{fonction} est réservé aux applications \emph{partielles} : une fonction $g:A\rightarrow B$ est définie sur son \emph{domaine de définition} $D_g\subset A$, sa restriction $g_{\vert D_g}$ à ce domaine étant une application $A\supset D_g \rightarrow B$. On peut aussi la voir comme une transition $g:A\rightsquigarrow B$ vérifiant : $\forall a\in A, card(g(a))\leq 1$. Comme pour les applications\footnote{Voir la remarque \ref{rmq appli det}.}, l'écriture $g(a)$ désignera en général dans ce cas, du moins si $a\in D_g$, un élément de $B$ (tandis que si $a\notin D_g$, on aura $g(a)=\emptyset$).
Les fonctions sont également appelées \emph{transitions quasi-déterministes}.
\end{rmq}

\subsubsection{Relations binaires, applications, fonctions}\label{subsubs relations binaires applications fonctions}

Une relation binaire $R$ d'un ensemble $S$ vers un ensemble $L$ est un triplet $(S,L,\vert R\vert)$, où $\vert R\vert\subset S\times L$ est le graphe de $R$, graphe que d'ailleurs nous nous autoriserons parfois à noter lui-même $R$. Pour tout $s\in S$, nous notons  
\[R(s)=\{l\in L, (s,l)\in \vert R\vert\},\]  et, pour tout $l\in L$,
\[R^{-1}(l)=\{s\in S, (s,l)\in \vert R\vert\}.\] L'\emph{image} $Im(R)$ est l'ensemble défini par $Im(R)=\bigcup_{s\in S} R(s)$, et le \emph{domaine de définition} $D_R$ est défini par $D_R=\{s\in S,R(s)\neq\emptyset\}$. L'ensemble des relations binaires de $S$ vers $L$ sera noté $\mathcal{B}_{(S,L)}$.

Conformément aux remarques \ref{rmq appli det} et \ref{rmq fonctions quasi det}, une \emph{fonction} $f$ de $S$ vers $L$ est une relation binaire $f$ de $S$ vers $L$ telle que pour tout $s\in S$ on a $card(R(s))\leq 1$, tandis qu'une \emph{application} $f$ de $S$ vers $L$ est une fonction dont le domaine de définition est $S$. 

\subsubsection{Relations multiples}\label{subsubs relations  multiples}

Une \emph{relation multiple}%
\footnote{Voir \cite{Dugowson:20150807}, section \textbf{§\,1.1}.} 
 $R$ d'index $I$ et de contexte $\mathcal{E}=(E_i)_{i\in I}$  est un triplet $(I,\mathcal{E},\vert R\vert)$ avec $\vert R\vert\subset \Pi_I\mathcal{E}$, le graphe de $R$, où l'on pose
\[
\Pi_I\mathcal{E}=\prod_{i\in I}{E_i}.
\]
Comme pour les relations binaires, on désignera parfois de la même façon une relation multiple et son graphe.\\
 
La classe des relations multiples de contexte $\mathcal{E}$ est notée $\mathcal{R}_\mathcal{E}$. Plus largement, la classe des relations multiples d'index $I$ est notée $\mathcal{R}_I$. En particulier, pour $I=2=\{0,1\}$, on retrouve les relations binaires.

\subsubsection{Relations binaires multiples}\label{subsubs relations binaires multiples}

Une \emph{relation binaire multiple}%
\footnote{Voir \cite{Dugowson:20150807}, section \textbf{§\,1.2}.} 
 $R$ d'index $I$, de contexte d'entrée $\mathcal{A}=(A_i)_{i\in I}$ et de contexte de sortie $\mathcal{B}=(B_i)_{i\in I}$ est un quadruplet  $(I,\mathcal{A},\mathcal{B}, \vert R\vert)$, avec $\vert R\vert\subset \Pi_I\mathcal{E}$ où $\mathcal{E}=(E_i=A_i\times B_i)_{i\in I}$ est parfois appelé le \emph{contexte global} de $R$. La classe des relations binaires multiples de contexte d'entrée $\mathcal{A}$ et de contexte de sortie $\mathcal{B}$ est notée  $\mathcal{BM}_{(\mathcal{A}, \mathcal{B})}$, et celle des relations binaires multiples \emph{non vides} de mêmes contextes sera notée $\mathcal{BM}^*_{(\mathcal{A}, \mathcal{B})}$. Plus largement, la classe de toutes les relations binaires multiples d'index $I$ est notée $\mathcal{BM}_I$.

Reprenons les opérateurs de \og transtypage\fg\, des relations binaires multiples $rd$, $rm$ et $rb$ définis   dans \cite{Dugowson:20150807}, section \textbf{§\,1.2.2}. Notant $2I=I\times \{0,1\}$, nous définissons un isomorphisme $rd:\mathcal{BM}_I  \simeq \mathcal{R}_{2I}$ en posant
\[
\mathcal{BM}_I\ni(I,\mathcal{A},\mathcal{B}, \vert R\vert)=R
\mapsto
rd(R)=(2I,\mathcal{D}, rd_\mathcal{E}(\vert R\vert))\in \mathcal{R}_{2I}, \]
avec
\begin{itemize}
\item $\mathcal{D}=(D_k)_{k\in 2I}$,
où pour tout $i\in I$,  $D_{(i,0)}=A_i$  et $D_{(i,1)}=B_i$,
\item $rd_\mathcal{E}(\vert R\vert)=\{rd_\mathcal{E}(u), u\in \vert R\vert\}$, avec, 
pour tout $u=\left((a_i,b_i)\right)_{i\in I}\in\Pi_I\mathcal{E} $, 
$rd_\mathcal{E}(u)=(d_k)_{k\in 2I}$, où  $d_{(i,0)}=a_i$ et $d_{(i,1)}=b_i$.\\
\end{itemize}

L'injection $rm:\mathcal{BM}_I\hookrightarrow \mathcal{R}_I$ est définie par
\[
\mathcal{BM}_I\ni(I,(A_i)_{i\in I},(B_i)_{i\in I}, \vert R\vert)=R
\mapsto
rm(R)=(I,(A_i\times B_i)_{i\in I}, \vert R\vert)\in \mathcal{R}_{I}, \]
et l'injection $rb:\mathcal{BM}_I\hookrightarrow \mathcal{R}_2$ par
\[
\mathcal{BM}_I\ni(I,\mathcal{A},\mathcal{B}, \vert R\vert)=R
\mapsto
rb(R)=(\Pi_I\mathcal{A}, \Pi_I\mathcal{B}, \vert rb(R)\vert)\in \mathcal{B}_{(\Pi_I\mathcal{A},\Pi_I\mathcal{B})}, \]
où le graphe $\vert rb(R)\vert$ est donné par un réarrangement évident des composantes des éléments de $\vert R\vert$.

\section{Dynamiques sous-fonctorielles ouvertes}\label{sec DSFO}

\subsection{Multi-dynamiques sous-fonctorielles}

\subsubsection{Définition des multi-dynamiques sous-fonctorielles}

On se donne une petite catégorie $\mathbf{C}$, et un ensemble non vide $M$.

\begin{df}\label{df multi-dyna sous-fonct}
Une \emph{multi-dynamique sous-fonctorielle} $\alpha$ de moteur $\mathbf{C}$ et d'ensemble paramétrique $M$  consiste en la donnée
\begin{itemize}
\item d'une application qui à tout objet $S\in\dot{\mathbf{C}}$ associe un ensemble $S^\alpha$, de telle sorte que $S\neq T \Rightarrow S^\alpha\cap T^\alpha=\emptyset$,
\item d'une application qui à toute flèche $(S\stackrel{f}{\rightarrow}T)\in\overrightarrow{\mathbf{C}}$ associe une famille de transitions $f^\alpha=(f^\alpha_\mu:S^\alpha\rightsquigarrow T^\alpha)_{\mu\in M}$ indexée par $M$ de telle sorte que pour tout $S\in\dot{\mathbf{C}}$ on ait%
\footnote{Sur la signification de l'expression $Id_{S^\alpha}$, voir la remarque \ref{rmq appli det} (page \pageref{rmq appli det}).}
\[
\forall \mu\in M, (Id_S)^\alpha_\mu\subset Id_{S^\alpha}
\]
et pour tout couple $(g,f)$ de flèches composables de $\mathbf{C}$, on ait
\[
\forall \mu\in M, (g\circ f)^\alpha_\mu\subset g^\alpha_\mu \odot f^\alpha_\mu.
\]
\end{itemize}

\end{df}

Nous écrirons $\alpha:\mathbf{C}\rightharpoondown \mathbf{P}^{\underrightarrow{\scriptstyle{M}}}$ pour indiquer que $\alpha$ est une telle multi-dynamique sous-fonctorielle.
L'écriture $f^\alpha: S^\alpha \twosquig_M  T^\alpha$ signifiera que $f^\alpha$ est la famille indexée par $M$ de transitions de $S^\alpha$ dans $T^\alpha$ associée à $f$ par $\alpha$. Nous pouvons ainsi écrire :

\begin{equation}\label{eq def multi-dynamique}
\begin{tabular}{rccc}
$\alpha :$  & $\mathbf{C}$ & $\rightharpoondown$ &  $\mathbf{P}^{\underrightarrow{\scriptstyle{M}}}$ \\ 
 & $(S\stackrel{f}{\rightarrow}T)$ & $ \longmapsto$ & $f^\alpha: S^\alpha \twosquig_M  T^\alpha$ \\ 
\end{tabular} 
\end{equation}

\paragraph{Ensemble des états.} Noté $st(\alpha)$, \emph{l'ensemble des états} de la dynamique $\alpha$ est défini par l'union (disjointe)
\[
st(\alpha)=\bigcup_{S\in \dot{\mathbf{C}}}S^\alpha.
\]

Pour tout état $s\in st(\alpha)$, nous noterons $typ(s)$ son type, autrement l'unique sommet $S\in\dot{\mathbf{C}}$ tel que $s\in S^\alpha$. Autrement dit, le type d'un état de la dynamique $\alpha$ est caractérisé par la relation
\[s\in (typ(s))^\alpha.\]

\paragraph{États \og hors-jeu\fg.}\label{etats hors jeu} Un état $s\in S^\alpha\subset st(\alpha)$ est dit \emph{hors-jeu pour le paramètre $\mu\in M$} si ${(Id_S)}_\mu^\alpha(s)=\emptyset$. Un état hors-jeu pour toutes les valeurs du paramètre sera bien entendu simplement dit \emph{hors-jeu}. Un état qui n'est pas hors-jeu sera dit \emph{dans le jeu}. Cette notion jouera un rôle crucial pour les \og dynamiques intemporelles\fg, c'est-à-dire celles de moteur $\mathbf{C}=\mathbf{1}$ (voir l'exemple \ref{exm dynamique intemporelle}).

\paragraph{Mono-dynamiques.} Si $M$ est réduit à un singleton, $\alpha$ est appelée une \emph{mono-dynamique sous-fonctorielle}, ou simplement  une \emph{dynamique (sous-fonctorielle)}.\\

\subsubsection{Multi-dynamiques graphiques}\label{subsubs multi dynamiques graphiques}

Étant donné un graphe
 $\mathbf{G}$, constitué d'un ensemble $\dot{\mathbf{G}}$ de sommets d'un ensemble $\overrightarrow{\mathbf{G}}$
d'arêtes dont chacune admet une source et un but parmi les sommets, et étant donné un ensemble non vide  $M$, on définit --- conformément à la définition 12 donnée dans \cite{Dugowson:20150807} --- une \emph{multi-dynamique graphique}  $\alpha$ de moteur $\mathbf{G}$ et d'ensemble paramétrique $M$  comme un morphisme de graphes%
\footnote{Pour la notation $\vert\mathbf{C}\vert$, voir les rappels de la section \textbf{§\,\ref{subsubs cat et graph}}.} %
  $\alpha:\mathbf{G}\longrightarrow \vert\mathbf{P}^{\underrightarrow{\scriptstyle{M}}}\vert$. Comme dans la définition \ref{df multi-dyna sous-fonct}, nous demanderons en outre que, pour tous sommets $S$ et $T$ de $\mathbf{G}$,  soit satisfaite la condition  $S\neq T \Rightarrow S^\alpha\cap T^\alpha=\emptyset$.
Autrement dit, une telle multi-dynamique graphique est 
constituée par la donnée
\begin{itemize}
\item d'une application qui à tout sommet $S\in\dot{\mathbf{G}}$ associe un ensemble $S^\alpha$ (avec $S\neq T \Rightarrow S^\alpha\cap T^\alpha=\emptyset$),
\item d'une application qui à toute arête $(S\stackrel{f}{\rightarrow}T)\in\overrightarrow{\mathbf{G}}$ associe une famille de transitions $(f^\alpha: S^\alpha \twosquig_M  T^\alpha)=(f^\alpha_\mu:S^\alpha\rightsquigarrow T^\alpha)_{\mu\in M}$ indexée par $M$.
\end{itemize}

\begin{rmq}\label{rmq condition graphique soit sous fonct} Soit $\mathbf{B}$ une petite catégorie, et $M$ un ensemble non vide.
\`{A} toute multi-dynamique sous-fonctorielle $\beta$ de moteur $\mathbf{B} $ on associe canoniquement, par oubli des propriétés sous-fonctorielles, une multi-dynamique graphique notée encore $\beta$ (ou, si l'on veut éviter toute ambiguïté, $\vert \beta\vert$) et de moteur le graphe
 $\vert\mathbf{B} \vert$.
Inversement, il est immédiat qu'une multi-dynamique graphique  \[\beta:\vert\mathbf{B} \vert\longrightarrow \vert\mathbf{P}^{\underrightarrow{\scriptstyle{M}}}\vert\] est sous-fonctorielle, et peut dès lors s'écrire 
\[\beta:\mathbf{B} \rightharpoondown\mathbf{P}^{\underrightarrow{\scriptstyle{M}}},\] si et seulement si elle vérifie les deux conditions suivantes :
\begin{itemize}
\item pour tout $S\in\dot{\mathbf{B} }$,
\[
 \forall \mu\in M, (Id_S)^\beta_\mu\subset Id_{S^\beta}
\]
\item pour tout couple $(g,f)$ de flèches composables de $\mathbf{B} $, 
\[
\forall \mu\in M, (g\circ f)^\beta_\mu\subset g^\beta_\mu \odot f^\beta_\mu.
\]
\end{itemize}
Nous aurions d'ailleurs pu définir de cette manière les multi-dynamiques sous-fonctorielles, ce qui aurait été en accord avec la façon dont les définitions \og sous-catégoriques\fg\, données dans \cite{Dugowson:20150809} s'appuient sur les définitions graphiques figurant dans \cite{Dugowson:20150807}.
\end{rmq}

\subsubsection{Multi-dynamorphismes}\label{subsubs multi dynamorphismes}

La catégorie des multi-dynamiques sous-fonctorielles a pour objets toutes les multi-dynamiques sous-fonctorielles $\alpha:\mathbf{C}\rightharpoondown \mathbf{P}^{\underrightarrow{\scriptstyle{L}}}$,
 où $\mathbf{C}$ décrit la classe des petites catégories et $L$ celle des ensembles, et pour flèches les multi-dynamorphismes définis de la façon suivante.

\begin{df}\label{df multi-dynamorphismes}
Étant données $\alpha:\mathbf{C}\rightharpoondown \mathbf{P}^{\underrightarrow{\scriptstyle{L}}}$ et
$\beta:\mathbf{D}\rightharpoondown \mathbf{P}^{\underrightarrow{\scriptstyle M}}$ deux multi-dynamiques, un \emph{multi-dynamorphisme} $(\theta,\Delta,\delta)$ de $\alpha$ vers $\beta$ consiste en la donnée
\begin{itemize}
\item d'une application $\theta:L\rightarrow M$,
\item d'un foncteur $\Delta:\mathbf{C}\rightarrow \mathbf{D}$,
\item d'une famille de transitions
$\delta=(\delta_S:S^\alpha\rightsquigarrow (\Delta S)^\beta)_{S\in \dot{\mathbf{C}}}$,
\end{itemize} 
\noindent tels que, pour tout $\lambda\in L$, $(\Delta,\delta)$ définit un \emph{mono-dynamorphisme} de $\alpha_\lambda$ vers $\beta_{\theta(\lambda)}$, ce qui signifie que pour tout $\lambda\in L$, tous $S$ et $T$ dans $\dot{\mathbf{C}}$ et tout $(e:S\rightarrow T)\in \overrightarrow{\mathbf{C}}$, on a
\[\delta_T\odot e^\alpha_\lambda\subset (\Delta e)^\beta_{\theta(\lambda)}\odot\delta_S.\] 
\end{df}

\begin{rmq}\label{rmq notation simplifiee des dynamorphismes}
En pratique, la famille $(\delta_S:S^\alpha\rightsquigarrow (\Delta S)^\beta)_{S\in \dot{\mathbf{C}}}$ sera identifiée à la transition $\delta:st(\alpha)\rightsquigarrow st(\beta)$ définie pour tout $s\in st(\alpha)$ par $\delta(s)=\delta_{typ(s)}(s)$, et la propriété ci-dessus reliant $\delta$, $\Delta$ et $\theta$ sera simplement écrite
\[\delta\odot e^\alpha_\lambda\subset (\Delta e)^\beta_{\theta(\lambda)}\odot\delta.\]
\end{rmq}

\paragraph{$\mathbf{C}$-multi-dynamorphismes.} Par convention, et sauf mention contraire, lorsque les deux dynamiques $\alpha$ et $\beta$ ont le même moteur $\mathbf{C}$, on appelle \emph{$\mathbf{C}$-multi-\-dyna\-morphisme}  de  $\alpha$ vers $\beta$ tout multi-dynamorphisme pour lequel, dans la définition ci-dessus, on a  $\Delta=Id_\mathbf{C}$.\\

En pratique, les multi-dynamorphismes seront simplement appelés des dynamorphismes.

\subsubsection{Quotient paramétrique}\label{subsubs quotient parametrique}

\begin{prop}\label{prop quotient param}
Soit $\alpha:\mathbf{C}\rightharpoondown\mathbf{P}^{\underrightarrow{\scriptstyle{M}}}$ une multi-dynamique sous-fonctorielle de moteur $\mathbf{C}$ et d'ensemble paramétrique $M$, et $\sim$ une relation d'équivalence sur $M$. Et soit $\beta$  la multi-dynamique graphique %
 de  moteur $\vert\mathbf{C}\vert$ et d'ensemble paramétrique $\widetilde{M}=M/{\sim}$ définie par
\begin{itemize}
\item $\forall S\in\vert\mathbf{C}\vert$,  $S^\beta=S^\alpha$,
\item $\forall(e:S\rightarrow T)\in\overrightarrow{\mathbf{C}}$, $\forall \lambda\in \widetilde{M}$, $\forall a\in S^\beta$, 
\[e^\beta_\lambda(a)=\bigcup_{\mu\in\lambda}e^\alpha_\mu(a).\]
\end{itemize}
Alors $\beta$ est une multi-dynamique sous-fonctorielle.
\end{prop}
\paragraph{Preuve.} Il suffit de vérifier que pour tout $\lambda\in \widetilde{M}$ la dynamique graphique $\beta_\lambda$ est sous-fonctorielle sur $\mathbf{C}$. Cela résulte du fait que $\beta_\lambda$ est l'union d'une famille de dynamiques sous-fonctorielles sur $\mathbf{C}$, une telle union étant sous-fonctorielle d'après la proposition 2, section 
\textbf{§\,2.1.7} de \cite{Dugowson:20150809}.
\begin{flushright}$\square$\end{flushright} 

\begin{df}\label{df reduc param multi}
La multi-dynamique sous-fonctorielle $\beta$ définie dans la proposition \ref{prop quotient param} est appelée \emph{quotient (paramétrique) de $\alpha$ par $\sim$} et on la note  $\alpha/{\sim}$.
\end{df}

\subsubsection{Multi-dynamiques (quasi-)déterministes} 

\begin{df} \label{df multi dyn deterministes} Soit $\alpha:\mathbf{C}\rightharpoondown \mathbf{P}^{\underrightarrow{\scriptstyle{M}}}$  une multi-dynamique sous-fonctorielle.
Si, pour \emph{toute} valeur paramétrique $\mu\in M$ et toute flèche $e\in\overrightarrow{\mathbf{C}}$, la transition $e_\mu^\alpha$ est déterministe%
\footnote{Autrement dit si $e_\mu^\alpha$ est une application; voir la remarque \ref{rmq appli det}.} %
(respectivement quasi-déterministe%
\footnote{Autrement dit si $e_\mu^\alpha$ est une fonction; voir la remarque \ref{rmq fonctions quasi det}.})
la multi-dynamique $\alpha$ est dite \emph{déterministe} (respectivement \emph{quasi-déterministe}).
 En outre, $\alpha$ est dite
 
\begin{itemize}
\item \emph{bien quasi-déterministe} si elle est quasi-déterministe mais n'est pas détermi\-niste,
\item \emph{pluraliste} si elle n'est pas quasi-déterministe.
\end{itemize} 
\end{df}

Bien entendu, toute multi-dynamique sous-fonctorielle déterministe est quasi-déterministe.

\subsubsection{Multi-dynamiques fonctorielles}

\begin{df}
Une \emph{multi-dynamique fonctorielle} $\alpha$ de moteur la petite catégorie $\mathbf{C}$ et d'ensemble paramétrique $M$ est un foncteur $\alpha:\mathbf{C}\rightarrow \mathbf{P}^{\underrightarrow{\scriptstyle{M}}}$, vérifiant en outre 
\[\forall (S,T)\in {\dot{\mathbf{C}}}^2, S\neq T \Rightarrow S^\alpha\cap T^\alpha=\emptyset.\] 
\end{df}

De façon équivalente, une multi-dynamique fonctorielle est une multi-dynamique sous-fonctorielle  $\alpha$ vérifiant en outre les deux relations
\[(Id_A)^\alpha=Id_{A^\alpha}\] et
\[
(f\circ e)^{\alpha} = f^{\alpha}\odot e^{\alpha},
\] pour tous choix de $A\in\dot{\mathbf{C}}$ et $(f,e)$ flèches composables de $\mathbf{C}$.

\begin{rmq}
La définition ci-dessus équivaut à celle d'une ``multi-dynamique catégorique propre'' donnée dans \cite{Dugowson:201112} et \cite{Dugowson:201203}. On fera toutefois attention au fait que le qualificatif \emph{propre} n'y a pas la même signification que dans les textes ultérieurs \cite{Dugowson:20150807} et \cite{Dugowson:20150809} : dans les premiers, il se rapporte aux dynamiques fonctorielles $\alpha$ pour lesquelles $S\neq T$ $ \Rightarrow$ $ S^\alpha\cap T^\alpha= \emptyset$, condition à présent automatiquement satisfaite, tandis qu'il se rapporte dans \cite{Dugowson:20150807} et \cite{Dugowson:20150809} aux dynamiques sous-fonctorielles  qui vérifient l'égalité $(Id_A)^\alpha=Id_{A^\alpha}$.
\end{rmq}

\paragraph{Dynamorphismes.} Ils sont définis comme dans le cas sous-fonctoriel. Aussi, les dynamorphismes entre deux multi-dynamiques fonctorielles sont les dynamorphismes entre les deux multi-dynamiques sous-fonctorielles sous-jacentes. Par conséquent, les premières constituent une sous-catégorie pleine de celle formée par les secondes. 

\subsubsection{Mono-dynamiques fonctorielles}

La notion de \og dynamique catégorique\fg\, développée dans \cite{Dugowson:201112} et \cite{Dugowson:201203}, où une dizaine d'exemples de ce type de dynamiques est donnée, coïncide avec celle de mono-dynamique fonctorielle. Ces dynamiques, non nécessai\-rement déterministes, constituent donc une sous-catégorie de celle des multi-dynamiques sous-foncto\-rielles.

\subsubsection{Horloges}

\begin{prop}\label{prop dysc deterministe implique fonctoriel}
Une dynamique sous-fonctorielle déterministe est nécessairement fonctorielle.
\end{prop}
\paragraph{Preuve.} Si $\alpha$ est déterministe, les deux membres des inclusions de la forme
\[(g\circ f)^{\alpha}(a) \subset (g^{\alpha}\odot f^{\alpha})(a)\] sont des singletons, il y a donc égalité.
\begin{flushright}$\square$\end{flushright}

\begin{df} Une \emph{horloge} $\mathbf{h}$ de moteur la petite catégorie $\mathbf{C}$ est un foncteur (covariant) \[ \mathbf{h}:\mathbf{C}\rightarrow \mathbf{Sets}\] satisfaisant en outre à la condition de distinction des états
\[\forall (S,T)\in \dot{\mathbf{C}}^2, 
S\neq T \Rightarrow S^\mathbf{h}\cap T^\mathbf{h}=\emptyset,\] 
où comme d'habitude nous posons
$S^\mathbf{h}=\mathbf{h}(S)$.
\end{df}

Autrement dit, une horloge  est une mono-dynamique fonctorielle déterministe sur $\mathbf{C}$.\\

\begin{exm}[Horloge essentielle, horloge existentielle]\label{exm horloge essentielle et existentielle}
 Dans \cite{Dugowson:201112} et \cite{Dugowson:201203} (sections \textbf{§\,3.4.3} et \textbf{§\,3.5.2}), nous associons fonctoriellement à toute petite catégorie $\mathbf{C}$ deux horloges particulières, appelées respectivement l'\emph{horloge essentielle} $\zeta_\mathbf{C}$ et l'\emph{horloge existentielle} $\xi_\mathbf{C}$ de $\mathbf{C}$. Rappelons que cette dernière est la $\mathbf{C}$-dynamique déterministe $\xi=\xi_\mathbf{C}$ telle que
 pour tout $T\in \dot{\mathbf{C}}$, on a
$T^\xi=\{\rightarrow T\}$,
et pour tout  $(f:S\rightarrow T)\in \overrightarrow{\mathbf{C}}$ et tout  $a\in S^\xi$, on a $f^\xi (a)=\{f\circ a\}$, où $\{\rightarrow S\}$ désigne la classe des flèches de but $S$ dans la catégorie considérée, $S$ désignant un objet de ladite catégorie.
\end{exm}

\begin{exm} Pour  $\mathbf{C}=\mathbf{R}_+$, et $t_0\in\mathbf{R}\cup\{-\infty\}$, l'ensemble des états $]t_0,+\infty[$ muni de l'action $\mathbf{h}$ de $\mathbf{R}_+$ définie pour tout $d\in\mathbf{R}_+$ et tout $t>t_0$ par $d^\mathbf{h}(t)=t+d$ constitue l'ensemble des instants d'une horloge. De même, pour $t_0\in\mathbf{R}$, l'ensemble des états $[t_0,+\infty[$ muni de l'action $\mathbf{h}$ définie comme ci-dessus est une horloge, l'horloge existentielle $\xi_{\mathbf{R}_+}$ correspondant au cas où $t_0=0$.
\end{exm}

\paragraph{Topos des horloges de moteur $\mathbf{C}$.}

En prenant pour flèches entre horloges de même moteur $\mathbf{C}$ les transformations naturelles, on définit une catégorie équivalente à $\mathbf{Sets}^\mathbf{C}$, topos des préfaisceaux d'ensembles sur $\mathbf{C}^{\mathrm{op}}$. En effet, on construit facilement une telle équivalence en associant canoniquement à tout foncteur $\mathbf{C}\rightarrow\mathbf{Sets}$ un foncteur équivalent mais satisfaisant la contrainte $\forall (S,T)\in \dot{\mathbf{C}}^2, S\neq T \Rightarrow S^\alpha\cap T^\alpha=\emptyset$.\\

Cela dit, les transformations naturelles entre horloges de moteur $\mathbf{C}$ sont des $\mathbf{C}$-dynamorphismes particuliers entre ces horloges vues comme multi-dynamiques, à savoir des 
\textit{dynamorphismes déterministes}, et 
il y a en général d'autres \emph{dynamorphismes} entre horloges que les seuls déterministes.

\paragraph{Catégorie des horloges.}\label{paragraph categorie des horloges}

On définit la \emph{catégorie des horloges}  en prenant pour objets toutes les horloges (pour tous les moteurs possibles), et pour flèches  entre deux horloges 
$\mathbf{h}:\mathbf{C}\rightarrow \mathbf{Sets}$ et $\mathbf{k}:\mathbf{D}\rightarrow\mathbf{Sets}$
tous les dynamorphismes \emph{quasi-déterministes}  $(\Delta,\delta):\mathbf{h}\looparrowright \mathbf{k}$. Autrement dit, conformément à la {remarque\,\ref{rmq fonctions quasi det}} (page \pageref{rmq fonctions quasi det}) et à  la {définition \ref{df multi-dynamorphismes}} (page \pageref{df multi-dynamorphismes}), une telle flèche consiste en 
un couple 
$(\Delta,\delta)$, avec pour 
$\Delta$ un foncteur $\mathbf{C}\rightarrow\mathbf{D}$ et pour $\delta$ une famille de transitions quasi-déterministes 
$(\delta_A: A^\mathbf{h}\rightsquigarrow A^\mathbf{k})_{A\in\dot{\mathbf{C}}}$  vérifiant la condition suivante
\[\forall (A\stackrel{e}{\rightarrow}B)\in \overrightarrow{\mathbf{C}},\, 
\delta_B\odot e^\mathbf{h} \subset (\Delta e)^\mathbf{k}\odot \delta_A.\] 

Conformément à la remarque \ref{rmq notation simplifiee des dynamorphismes} (page \pageref{rmq notation simplifiee des dynamorphismes}), 
la transition $\delta_A$ sera simplement notée $\delta$, de sorte que, pour toute flèche $e\in\overrightarrow{\mathbf{C}}$, 
la relation que doivent satisfaire ces transitions s'écrit : $\delta\odot e^\mathbf{h} \subset (\Delta e)^\mathbf{k}\odot \delta$.

Dans le cas d'un dynamorphisme déterministe, autrement dit lorsque $\delta$ est une application, cette dernière condition s'écrit plus simplement 
\begin{equation}\label{eq dynamorphisme deterministe entre horloges}
\delta\circ e^\mathbf{h} = (\Delta e)^\mathbf{k}\circ \delta.
\end{equation}

Par ailleurs, lorsque $\mathbf{D}=\mathbf{C}$, on ajoute à la définition de la catégorie des horloges la condition $\Delta=Id_{\mathbf{C}}$ pour définir les flèches de la \emph{catégorie des horloges de  moteur $\mathbf{C}$}. En se restreignant de plus aux dynamorphismes déterministes, nous retrouvons les flèches du topos des horloges de moteur $\mathbf{C}$.

\paragraph{Instants et antériorité.}\label{subs instants et anteriorite}

Les états d'une horloge $\mathbf{h}$ sont appelés ses \emph{instants}. Une relation de pré-ordre, appelée \emph{antériorité}, est ainsi définie entre les instants d'une horloge $\mathbf{h}$ : $s$ est antérieur à $t$, ce que l'on note $s\leq_\mathbf{h} t$, si et seulement si
\[\exists e\in \overrightarrow{\mathbf{C}}, e^\mathbf{h}(s)=t.\] 

\subsection{Dynamiques sous-fonctorielles ouvertes}

\subsubsection{Définition}

\begin{df}\label{df DySCO}  Une \emph{dynamique sous-fonctorielle ouverte} $A$ de moteur $\mathbf{C}$ est la donnée 
\[A=\left((\alpha:\mathbf{C}\rightharpoondown \mathbf{P}^{\underrightarrow{\scriptstyle{M}}}) \stackrel{\tau}{\looparrowright}  (\mathbf{h}:\mathbf{C}\rightarrow \mathbf{P})\right)\]
\begin{itemize}
\item d'un ensemble non vide $M$,
\item d'une $\mathbf{C}$-multi-dynamique sous-fonctorielle $\alpha:\mathbf{C}\rightharpoondown \mathbf{P}^{\underrightarrow{\scriptstyle{M}}}$,
\item d'une $\mathbf{C}$-horloge $\mathbf{h}$,
\item d'un $\mathbf{C}$-multi-dynamorphisme \emph{déterministe}
\[
\tau:
(\alpha:\mathbf{C}\rightharpoondown \mathbf{P}^{\underrightarrow{\scriptstyle{M}}} )
\looparrowright 
(\mathbf{h}:\mathbf{C}\rightarrow \mathbf{P}).\]
\end{itemize}

La dynamique $A$ sera parfois désignée par sa partie multi-dynamique $\alpha$. En particulier, les états de $A$ sont définis comme ceux de $\alpha$, et l'on écrira $st(A)=st(\alpha)$. Le dynamorphisme $\tau$ est appelé la \emph{scansion} ou la \emph{datation} de $A$.

Une telle dynamique sous-fonctorielle ouverte est dite \emph{fonctorielle} si $\alpha$ est une multi-dynamique fonctorielle.
\end{df}

\subsubsection{Dynamorphismes entre dynamiques sous-fonctorielles ouvertes}

Conformément à \cite{Dugowson:20150809}, section \textbf{§\,2.4.2}, nous définissons ainsi les dynamorphismes entre dynamiques sous-fonctorielles ouvertes :

\begin{df}\label{df dyna DySCO} \label{df dyna DynO} On appelle \emph{dynamorphisme} d'une dynamique sous-fonctorielle ouverte
\[
A=
(
\rho:
(\alpha:\mathbf{C}\rightharpoondown \mathbf{P}^{\underrightarrow{\scriptstyle L}})
\looparrowright 
(\mathbf{h}:\mathbf{C}\rightarrow \mathbf{P})
)
\]
vers une dynamique sous-fonctorielle ouverte
\[
B=
(
\tau:
(\beta:\mathbf{D}\rightharpoondown \mathbf{P}^{\underrightarrow{\scriptstyle M}})
\looparrowright 
(\mathbf{k}:\mathbf{D}\rightarrow \mathbf{P})
)
\]
la donnée d'un quadruplet $(\theta,\Delta,\delta,d)$ tel que
\begin{enumerate}
\item $(\theta,\Delta,\delta)$ est un multi-dynamorphisme sous-fonctoriel de $\alpha$ vers $\beta$,
\item $(\Delta,d)$ est un dynamorphisme de $\mathbf{h}$ vers $\mathbf{k}$,
\item pour tout $S\in\dot{\mathbf{C}}$, la condition suivante de synchronisation entre $\rho$ et $\tau$ est satisfaite :
\[\tau_{\Delta_S}\odot \delta_S\subset d_S\odot \rho_S. \]  
\end{enumerate}

\end{df}

Conformément à la définition 16 de \cite{Dugowson:20150807}, étant donné un dynamorphisme $(\theta,\Delta,\delta,d)$, on appellera
\begin{itemize}
\item $\theta$ sa  \emph{partie paramétrique},
\item $\Delta$ sa  \emph{partie fonctorielle},
\item $\delta$ sa  \emph{partie transitionnelle},
\item et $d$ sa \emph{partie horloge}.
\end{itemize}

En prenant pour flèches les dynamorphismes, la classe des dynamiques sous-fonctorielles ouvertes constitue une catégorie notée $\mathbf{DySCO}$ dans \cite{Dugowson:20150809}.

\subsection{Réalisations d'une dynamique sous-fonctorielle ouverte}

\subsubsection{Définitions}

\begin{df}
\'{E}tant donnée
\[
A=(
\tau:
(\alpha=(\alpha_\lambda)_{\lambda \in L}:\mathbf{C}\rightharpoondown \mathbf{P}^{\underrightarrow{\scriptstyle{L}}} )
\looparrowright 
(\mathbf{h}:\mathbf{C}\rightarrow \mathbf{P})
)
\]
une $\mathbf{C}$-dynamique ouverte, une
\emph{réalisation} $\mathfrak{a}$ de $A$ consiste en un $\mathbf{C}$-dynamorphisme quasi-déterministe $\mathfrak{a}:\mathbf{h}\looparrowright A$ tel que $\tau\odot\mathfrak{a}\subset Id_\mathbf{h}$.
\end{df}

On vérifie facilement%
\footnote{Voir les section \textbf{§\,2.5} de \cite{Dugowson:20150807} et \cite{Dugowson:20150809}}
 qu'une telle réalisation $\mathfrak{s}$ de $A$ consiste en
un couple $\mathfrak{s}=(\lambda,\mathfrak{a})$ constitué d'une valeur $\lambda\in L$ et d'une fonction $\mathfrak{a}:st(\mathbf{h})\dashrightarrow st(\alpha)$ définie sur une partie $D_\mathfrak{a}\subset st(\mathbf{h})$ et qui vérifie les  propriétés suivantes :

\[\forall t\in D_\mathfrak{a}, \tau(\mathfrak{a}(t))=t,\]
\[\forall S\in\dot{\mathbf{C}},\forall t\in S^\mathbf{h}\cap  D_\mathfrak{a}, \mathfrak{a}(t)\in S^\alpha,\]
\[\forall (S\stackrel{f}{\rightarrow}T)\in\overrightarrow{\mathbf{C}}, 
 \forall t\in S^\mathbf{h}, 
\left( f^\mathbf{h}(t)\in D_\mathfrak{a}\Rightarrow 
\left( t\in D_\mathfrak{a} 
\,\mathrm{et}\,
\mathfrak{a}(f^\mathbf{h}(t))\in f^\alpha_\lambda (\mathfrak{a}(t))\right)\right).\]

\begin{rmq}\label{rmq domaine def des real} La caractérisation donnée ci-dessus explicite en particulier le fait que si une réalisation $\mathfrak{a}$ de $A$ est définie à un instant $s$, alors elle est également définie à tout instant antérieur $t\leq s$, puisqu'il existe alors $f\in\overrightarrow{\mathbf{C}}$ tel que $s=f^\mathbf{h}(t)$. 
\end{rmq}

Le paramètre $\lambda$ est appelé la \emph{partie interne} ou \emph{paramétrique} de la réalisation $(\lambda,\mathfrak{a})$ de $A$, tandis que $\mathfrak{a}$ est la \emph{partie externe} de cette réalisation. Dans le cas où $A$ est une mono-dynamique ouverte, autrement dit si $L$ est un singleton, une réalisation  s'identifie à sa partie externe puisque la partie paramétrique est nécessairement égale à l'unique élément de $L$.

Nous notons $\mathcal{S}_A$ l'\emph{ensemble des parties externes des réalisations} de la dynamique sous-fonctorielle ouverte $A$, et $\mathcal{S}_{(A,\lambda)}$ ou $\mathcal{S}_{A_\lambda}$ l'ensemble des réalisations de la mono-dynamique ouverte $A_\lambda=(
\tau:
(\alpha_\lambda:\mathbf{C}\rightharpoondown \mathbf{P})
\looparrowright 
(\mathbf{h}:\mathbf{C}\rightarrow \mathbf{P})
)$, de sorte que 
\[\mathcal{S}_A=\bigcup_{\lambda\in L}\mathcal{S}_{A_\lambda},\]
cette union n'étant pas en général disjointe.

\begin{rmq}\label{rmq usage realisation au lieu de partie externe} Souvent, et sans que  cela ne porte  à conséquence, nous parlerons de \emph{la réalisation $\mathfrak{a}$ de $A$} au lieu de 
\emph{la partie externe $\mathfrak{a}$ d'une réalisation de $A$}. Cette façon de parler conduira par exemple à désigner l'ensemble $\mathcal{S}_A$ comme \og ensemble des réalisations de $A$\fg\, bien que l'expression soit impropre et qu'il vaille parfois mieux l'éviter.
\end{rmq}

\begin{rmq}
Quelle que soit la dynamique $A$, on a  $\mathcal{S}_A\neq\emptyset$, comme le prouve l'exemple ci-dessous.
\begin{exm} On appelle \emph{réalisation vide} de $A$ toute réalisation de $A$ dont la partie externe est vide, autrement dit tout couple de la forme $(\lambda,\underline{\emptyset})$, 
où $\lambda\in L$ et $ \underline{\emptyset}$ désigne la \emph{fonction} vide $st(H)\supset D_{\underline{\emptyset}}=\emptyset\hookrightarrow st(A)$. Quelle que soit la dynamique $A$, l'ensemble de ses réalisations vides est non vide, isomorphe à $L$. Toutes les réalisations vides de $A$ ont la même partie externe, notée $\underline{\emptyset}_A$ ou simplement 
$\underline{\emptyset}$,
que nous appellerons \emph{la réalisation vide} de $A$. On a donc toujours  $\mathcal{S}_A\ni\underline{\emptyset}_A$.
\end{exm} 

\paragraph{Notation.} L'ensemble des réalisations non vides de $A$ sera noté $\mathcal{S}_A^*$ :
\[\mathcal{S}_A^*=\mathcal{S}_A\setminus \{\underline{\emptyset}_A\}.\]
\end{rmq}

\begin{df}\label{df dynamique efficiente} Une dynamique sous-fonctorielle ouverte $A$ n'admettant pour seule réalisation que la réalisation vide, autrement dit telle que $\mathcal{S}_A^*=\emptyset$, sera dite \emph{inefficiente}. Dans le cas contraire, elle sera dite efficiente.
\end{df}

\subsubsection{Réalisations passant par un état}
\label{subsubs realisations sc passant par etat}

\begin{df}
Étant donnée une dynamique sous-fonctorielle ouverte $A$,  nous dirons qu'une réalisation\footnote{Voir la remarque \ref{rmq usage realisation au lieu de partie externe} ci-dessus.} $\mathfrak{a}$ de $A$ passe par un état $a\in st(A)$, si  $\mathfrak{a}(\tau(a))=a$.
\end{df}

Nous écrirons
\[\mathfrak{a}\rhd a,\] 
pour exprimer que $\mathfrak{a}$ passe par $a$.
 
 Ainsi, pour $A=(\tau:(\alpha_\lambda)_{\lambda\in L}\looparrowright \mathbf{h})$, on a 
\[\mathfrak{a}\rhd a \Leftrightarrow \mathfrak{a}(\tau(a))=a.\]

Plus généralement, si $E$ est un ensemble d'états de la dynamique ouverte $A$, nous écrirons 
\[\mathfrak{a}\rhd E\]
pour exprimer que $\mathfrak{a}$ passe par chacun des états $a\in E$. Dans le cas où $E$ est un ensemble fini $E=\{a_1,..., a_n\}$, nous écrirons souvent
\[\mathfrak{a}\rhd a_1,..., a_n\] au lieu de $\mathfrak{a}\rhd E$.

\begin{rmq}\label{rmq notation a passe par a et b} Par rapport à \cite{Dugowson:20150807} et \cite{Dugowson:20150809}, nous avons  légèrement changé la signification de l'expression $\mathfrak{a}\rhd a, b$. En effet, nous avons à présent
\[ (\mathfrak{a}\rhd a, b) \Leftrightarrow (\mathfrak{a}\rhd a \,\mathrm{et}\, \mathfrak{a}\rhd  b),\] alors que dans \cite{Dugowson:20150807} et \cite{Dugowson:20150809} l'écriture $\mathfrak{a}\rhd a, b$ signifiait non seulement que $\mathfrak{a}$ passait par $a$ et passait par $b$, mais aussi que $\tau(a)$ était antérieur\footnote{Voir plus haut la section \textbf{§\,\ref{subs instants et anteriorite}}.}  à $\tau(b)$. Désormais, cette dernière condition n'est donc plus  requise.
\end{rmq}

\subsection{Une classification des dynamiques sous-fonctorielles ouvertes}

Soit 
\[A=\left((\alpha:\mathbf{C}\rightharpoondown \mathbf{P}^{\underrightarrow{\scriptstyle{M}}}) \stackrel{\tau}{\looparrowright}  (\mathbf{h}:\mathbf{C}\rightarrow \mathbf{P})\right)\]
une dynamique sous-fonctorielle ouverte.

\paragraph{Type paramétrique.}

Nous dirons que $A$ est 
\begin{itemize}
\item \emph{paramétrique}, ou de type $\overline{\pi}$, si $card(M)>1$, 
\item \emph{non paramétrique}, ou de type $\dot{\pi}$, si $card(M)=1$.
\end{itemize}

Notons que $A$ est nécessairement soit de type $\dot{\pi}$, soit de type  $\overline{\pi}$.

\paragraph{Type de déterminisme.}

$A$ sera dite \emph{déterministe} (resp. \emph{quasi-déterministe}, \emph{bien quasi-déterministe}, \emph{pluraliste}) si $\alpha$ l'est.\\

Nous dirons en outre que $A$ est 
\begin{itemize}
\item de type $\overline{\delta}$ si elle est pluraliste.
\item de type $\delta$ si elle est déterministe,
\item de type $\underaccent{\dot}{\delta}$, si elle est bien quasi-déterministe.
\end{itemize}

Notons que $A$ est nécessairement soit de type $\overline{\delta}$, soit de type $\delta$, soit de type $\underaccent{\dot}{\delta}$.

\paragraph{Type de fonctorialité.}

Nous dirons que $A$ est 
\begin{itemize}
\item \emph{bien sous-fonctorielle}, ou de type $\underline{\phi}$, si $\alpha$ est sous-fonctorielle mais non fonctorielle,
\item \emph{fonctorielle non déterministe}, ou de type ${\phi}$, si $\alpha$ est fonctorielle mais non déterministe.
\end{itemize}

\begin{rmq}
Une dynamique ouverte de type $\underline{\phi}$ est nécessairement non déterministe d'après la proposition \ref{prop dysc deterministe implique fonctoriel}. Ainsi, $A$ est nécessairement soit de type $\underline{\phi}$, soit de type $\delta$, soit de type $\phi$. 
\end{rmq}

\begin{df} 
Nous dirons que la dynamique sous-fonctorielle ouverte $A$ est de type $[PDF\mathbf{C}]$, où les lettres $P$, $D$ et $F$ représentent des symboles pris respectivement dans les ensembles suivants 
\begin{itemize}
\item  $P\in\{\dot{\pi},\overline{\pi}\}$,
\item  $D\in\{\underaccent{\dot}{\delta}, \delta,\overline{\delta}\}$, 
\item $F\in\{\phi,\underline{\phi}\}$,
\end{itemize}
pour exprimer que $\mathbf{C}$ est le moteur de $A$,  $P$ est son type paramétrique, $D$ est son type de déterminisme et, dans le cas où $D\neq\delta$, $F$ est son type de fonctorialité. Si $D=\delta$, la place de $F$ est laissée vide.
\end{df}

La partie $[PDF]$ de cette classification permet de distinguer dix types de dynamiques sous-fonctorielles ouvertes, à savoir cinq types de dynamiques non paramétriques  

{\setlength{\baselineskip}{1.35\baselineskip}
\begin{itemize}
\item $[\dot{\pi}\underaccent{\dot}{\delta} \phi]$  : fonctorielles bien quasi-déterministes,
\item $[\dot{\pi}\underaccent{\dot}{\delta}\underline{\phi}]$  : bien sous-fonctorielles et bien quasi-déterministes,
\item $[\dot{\pi}\delta]$  : déterministes,
\item $[\dot{\pi}\overline{\delta}{\phi}]$  : fonctorielles pluralistes,
\item $[\dot{\pi}\overline{\delta}\underline{\phi}]$  : bien sous-fonctorielles et pluralistes,
\end{itemize}
\par}
et cinq types de dynamiques paramétriques
{\setlength{\baselineskip}{1.35\baselineskip}
\begin{itemize}
\item $[\overline{\pi}\underaccent{\dot}{\delta} \phi]$  : fonctorielles et bien quasi-déterministes,
\item $[\overline{\pi}\underaccent{\dot}{\delta}\underline{\phi}]$  : bien sous-fonctorielles et bien quasi-déterministes,
\item $[\overline{\pi}\delta]$  : déterministes,
\item $[\overline{\pi}\overline{\delta}{\phi}]$   : fonctorielles pluralistes,
\item $[\overline{\pi}\overline{\delta}\underline{\phi}]$  : bien sous-fonctorielles et pluralistes.
\end{itemize}
\par}

\begin{rmq} Si le moteur est un groupe (ou un groupoïde) $\mathbf{G}$ et que la dynamique considérée est fonctorielle, celle-ci est de façon évidente nécessairement déterministe. Par conséquent, les dynamiques de moteur $\mathbf{G}$ sont soit détermi\-nistes, 
soit de type 
$[\dot{\pi}\underaccent{\dot}{\delta}\underline{\phi}\mathbf{G}]$, 
$[\dot{\pi}\overline{\delta}\underline{\phi}\mathbf{G}]$, 
$[\overline{\pi}\underaccent{\dot}{\delta}\underline{\phi}\mathbf{G}]$ 
ou  
$[\overline{\pi}\overline{\delta}\underline{\phi}\mathbf{G}]$ (mais jamais de type $[\dot{\pi}\underaccent{\dot}{\delta} \phi\mathbf{G}]$, $[\dot{\pi}\overline{\delta} \phi\mathbf{G}]$, $[\overline{\pi}\underaccent{\dot}{\delta} \phi\mathbf{G}]$ ou $[\overline{\pi}\overline{\delta} \phi\mathbf{G}]$).
En outre, dans le cas où $\mathbf{G}=\mathbf{1}=\{\bullet\}$, la seule flèche de $\mathbf{1}$ étant une identité, les dynamiques de moteur $\mathbf{1}$ ne peuvent être pluralistes et sont donc soit déterministes, soit de type $[\dot{\pi}\underaccent{\dot}{\delta}\underline{\phi}\mathbf{1}]$, soit de type $[\overline{\pi}\underaccent{\dot}{\delta}\underline{\phi}\mathbf{1}]$.
\end{rmq}

\subsection{Quelques exemples}


\subsubsection{Une \og source \fg\, ($\mathbb{Y}$ = \textcjheb{y})}\label{subsubs exm futur}

\begin{exm}[{$\mathbb{Y}$ : une \og source lipschitzienne\fg}]\label{exm source lip}

Nous appelons \og source lipschitzienne\fg, la dynamique 
\[
{\mathbb{Y}}=
(
{\tau_{\mathbb{Y}}}:
({\alpha_{\mathbb{Y}}}:{\mathbf{C}_{\mathbb{Y}}}\rightharpoondown \mathbf{P}^{\underrightarrow{\scriptstyle {L_{\mathbb{Y}}}}})
\looparrowright 
({\mathbf{h}_{\mathbb{Y}}}:{\mathbf{C}_{\mathbb{Y}}}\rightarrow \mathbf{P})
),
\] 
dont
\begin{itemize}
\item le moteur est $\mathbf{C}_{\mathbb{Y}}=\mathbf{R}_+$, 
\item l'horloge est l'horloge existentielle $\mathbf{h}_{\mathbb{Y}}=\xi_{\mathbf{R}_+}$, \emph{i.e.} : $\bullet^{\mathbf{h}_{\mathbb{Y}}}=\mathbf{R}_+$ et $d^{\mathbf{h}_{\mathbb{Y}}}(t)=t+d$,
\item la paramétrisation est $L_{\mathbb{Y}}=\{*\}$,
\item l'ensemble des états est $st({\mathbb{Y}})=st(\alpha_{\mathbb{Y}})=\mathbf{R}_+\times \mathbf{R}$,
\item la scansion est donnée par $\tau_{\mathbb{Y}}(t,a)=t$,
\item et la loi s'écrit $d^{\mathbb{Y}}(t,a)=d^{\alpha_{\mathbb{Y}}}(t,a)=\{t+d\}\times [a-d,a+d]$.
\end{itemize}

\mbox{}

On vérifie aisément que l'ensemble des réalisations de $\mathbb{Y}$ est 
\[
\mathcal{S}_{\mathbb{Y}}=
Lip^1_+,
\]
et, puisque $\underline{\emptyset}_\mathbb{Y}=\emptyset$, l'ensemble des réalisations non vides de $\mathbb{Y}$ est
\begin{equation}\label{eq ens real non vide de Y}
\mathcal{S}_{\mathbb{Y}}^*=
(Lip^1_+)^*.
\end{equation}

C'est une dynamique de type  $[{\dot{\pi}}\overline{\delta}\phi\mathbf{R}_+]$, autrement dit fonctorielle, non paramétrique, pluraliste, et de moteur $\mathbf{R}_+$.

\begin{rmq}
Dans notre travail \cite{Dugowson:20170101}  inspiré par  la théorie \emph{métachronologique} du philosophe Pierre Michel Klein \cite{KleinPM:2014}, la dynamique $\mathbb{Y}$ est notée \textcjheb{y} --- \emph{yod}, en hébreu\footnote{ En \TeX, nous codons la lettre \textcjheb{y}  avec un \texttt{y} (code \texttt{\textbackslash textcjheb\{y\}}). } ---  du fait qu'elle y joue un rôle moteur aux côtés de la dynamique \og intemporelle\fg\, que nous présentons plus loin (exemple \ref{exm dynamique intemporelle}) et qui est également notée par une lettre hébraïque --- \textcjheb{w} (\emph{vav}) --- en référence au rôle essentiel joué par cette lettre dans la théorie de P. M. Klein.
\end{rmq}

\end{exm}


\subsubsection{Une \og histoire\fg\, ($\mathbb{H}$ = \textcjheb{h})}\label{subsubs exm passe}

\begin{exm}[$\mathbb{H}$, une dynamique \og historique \fg]\label{exm dyna historique}
 \'{E}tant donnés
\begin{itemize}
\item un indice $k\in\mathbf{N}\cup \{+\infty\}$,
\item un élément $a_0\in\{-\infty\}\cup\mathbf{R}$, appelé \emph{origine des histoires}, 
\item un élément $T_0\in\{-\infty\}\cup\mathbf{R}$, appelé \emph{origine des temps}, 
\item un ensemble non vide $L\subset \mathcal{P}(]T_0,+\infty[\times \mathbf{R})$ de parties du (demi-)plan $t>T_0$,
\end{itemize}
on 
appelle \emph{\og dynamique historique\fg\, de classe $\mathcal{C}^k$, d'origine des histoires $a_0$, d'origine des temps $T_0$ et de paramétrage $L$} la dynamique notée $\mathbb{H}_{(k,a_0, T_0, L)}$, ou simplement $\mathbb{H}$ s'il n'y a pas d'ambiguïté, 
 définie par
\[
{\mathbb{H}}=
(
{\tau_{\mathbb{H}}}:
({\alpha_{\mathbb{H}}}:{\mathbf{C}_{\mathbb{H}}}\rightharpoondown \mathbf{P}^{\underrightarrow{\scriptstyle{L_{\mathbb{H}}}}})
\looparrowright 
({\mathbf{h}_{\mathbb{H}}}:{\mathbf{C}_{\mathbb{H}}}\rightarrow \mathbf{P})
),
\] 
avec 
\begin{itemize}
\item pour moteur : $\mathbf{C}_{\mathbb{H}}=\mathbf{R}_+$, 
\item pour horloge : $\bullet^{\mathbf{h}_{\mathbb{H}}}=]T_0,+\infty[$ et, pour tout $d\in\mathbf{R}_+$ et tout $t>T_0$, $d^{\mathbf{h}_{\mathbb{H}}}(t)=t+d$,
\item pour paramétrisation : 
$L_{\mathbb{H}}=L$,
\item pour ensemble d'états%
\footnote{\`{A} propos de la notation $\mathcal{C}^k(]a,b[)$, voir la remarque \ref{rmq Ck(a,b) avec a>b}, section \textbf{§\,\ref{subsubs reels}}.}
 : $st({\mathbb{H}})=st(\alpha_{\mathbb{H}})=
\bigcup_{t\in]T_0,+\infty[}\left(\{t\}\times\mathcal{C}^k(]a_0,t[)\right)$, autrement dit 
\[
st({\mathbb{H}})=
\{(t,f), t>T_0\mathrm{\,et\,} f\in \mathcal{C}^k(]a_0,t[)\},
\]
\item pour scansion : $\tau_{\mathbb{H}}(t,f)=t$,
\item et pour loi celle définie pour tout $\gamma\in L$, tout $(t,f)\in st(\alpha_{\mathbb{H}})$ et tout $d\in\mathbf{R}_+$ par
\begin{equation}\label{eq transitions historiques}
d^{\mathbb{H}}_\gamma(t,f)=
\left\lbrace
(t+d,g)\in st({\mathbb{H}}), 
\left\vert
\begin{tabular}{l}
$g_{\vert ]a_0,t[}=f$, \\ 
$\{(s,g(s)), s\in ]t,t+d[\cap D_g\}\subset \gamma$.\\ 
\end{tabular} 
\right.
\right\rbrace
\end{equation}
\end{itemize}

Le type de la  dynamique $\mathbb{H}$ dépend du choix des réels $a_0$ et $T_0$ et de l'ensemble $L$. 

\paragraph{Considérons le cas où $k=1$,
 $a_0=-\infty$, $T_0=0$  et $L={\mathcal{GC}_+^*}$.}

Dans ce cas, où ${\mathcal{GC}_+^*}$ désigne l'ensemble des graphes des fonctions numériques continues dont le domaine de définition est un intervalle ouvert inclus dans $\mathbf{R}_+^*$ et admettant $0$ pour borne inférieure s'il n'est pas vide, autrement dit 
\[{\mathcal{GC}_+^*}=
\left\lbrace
\gamma=\{(s,l_\gamma(s)), s\in ]0,r_\gamma[\}, 
r_\gamma\in\overline{\mathbf{R}}_+, 
l_\gamma\in\mathcal{C}(]0,r_\gamma[)
\right\rbrace,\]
nous obtenons une dynamique de type 
$[\overline{\pi}\underaccent{\dot}{\delta}{\phi}\mathbf{R}_+]$.
Dans la suite, cette dynamique $\mathbb{H}_{(k=1,a_0-\infty, T_0=0, L={\mathcal{GC}_+^*})}$ sera simplement notée $\mathbb{H}$. 

\subparagraph{Identification des graphes aux fonctions.} On a $\gamma\in {\mathcal{GC}_+^*}$ si $\gamma$ est le graphe d'une fonction numérique $l$ 
définie et continue sur un intervalle de la forme $]0,r[$ avec $r\in \overline{\mathbf{R}}_+$.
 En identifiant le paramètre $\gamma$ à l'unique fonction continue $g_\gamma$ dont il est le graphe on peut écrire
\[{\mathcal{GC}_+^*}\simeq
\bigcup_{r\in \overline{\mathbf{R}}_+}{\mathcal{C}(]0,r[)},\]
le cas où $\gamma=\emptyset$ correspondant à $r=0$.\\

Décrivons les réalisations de $\mathbb{H}$. Pour cela, donnons-nous d'abord un $\gamma\in L$, posons pour simplifier $g=g_\gamma$ et notons $]0,r[=D_g$ son domaine de définition, où $r\in \overline{\mathbf{R}}_+$. 
Une réalisation non vide%
\footnote{On remarque que si $\gamma=\emptyset$, autrement dit si $r=0$,  $\mathbb{H}_\gamma$  n'admet pas de réalisation non vide. En effet, pour toute réalisation $\mathfrak{h}$ qui serait non vide, il existe $t>0$ tel que $\mathfrak{h}(t)$ soit défini, ce qui implique --- d'après la remarque \ref{rmq domaine def des real} --- que pour tout $s\in]0,t[$,  $\mathfrak{h}(s)$ est également défini. On doit alors avoir $\mathfrak{h}(t)\in (t-s)^\mathbb{H}_\gamma(\mathfrak{h}(s))$, mais d'après la formule (\ref{eq transitions historiques}) ceci est impossible lorsque $\gamma=\emptyset$. De même, $\mathbb{H}_\gamma$ n'admet pas de réalisation non vide si $g$ n'est de classe $C^1$ sur aucun intervalle de la forme $]0,t[$.}
 $\mathfrak{h}$ de $\mathbb{H}_\gamma$ 
est caractérisée par la donnée d'un couple $(D,h)$ constitué
\begin{itemize}
\item du domaine de définition $D=D_{\mathfrak{h}}$ de $\mathfrak{h}$, qui est un intervalle non vide  $D\subset \mathbf{R}_+^*$  de la forme 
$]0,a[$ ou\footnote{Si $a=r=+\infty$, il n'y a qu'une seule forme possible pour $D$, à savoir $D= \mathbf{R}_+^*$.} %
 $]0,a]$, avec $a\in ]0,r]$,
\item d'une application  $h:]-\infty,a[\rightarrow \mathbf{R}$ de classe $\mathcal{C}^1$  telle que $h_{\vert ]0,a[}=g_{\vert ]0,a[}$,
\end{itemize}
la réalisation  $\mathfrak{h}$ caractérisée par un tel couple $(D,h)$ étant alors  l'application  $D\rightarrow st(\mathbb{H})$ définie par 
\begin{equation}\label{eq solution (D,f) de Histoire}
\forall t\in D, \mathfrak{h}(t)=(t,h_{\vert ]-\infty,t[}).
\end{equation}

Par conséquent, à toute (partie externe d'une) réalisation non vide  $\mathfrak{h}\in\mathcal{S}_\mathbb{H}^*$ se trouve associé un tel couple $(D,h)$, unique, avec $h\in\mathcal{C}^1(\mathbf{R})$. Nous noterons en particulier $\mathfrak{h}\mapsto \widetilde{\mathfrak{h}}=h$ l'application de $\mathcal{S}_\mathbb{H}^*$ dans $\mathcal{C}^1(\mathbf{R})$ qui à une telle réalisation $\mathfrak{h}$ associe la fonction $h$ correspondante.


Mise à part la réalisation vide $\underline{\emptyset}_\mathbb{H}=\emptyset:\emptyset\hookrightarrow st(\mathbb{H})$, nous avons donc deux sortes de réalisations, selon que l'intervalle $D$ est ouvert ou fermé en sa borne supérieure. Celles pour lesquelles $D$ est ouvert forment un ensemble qui s'identifie à\footnote{
Voir les notations introduites en  section \textbf{§\,\ref{subsubs notations classes applications}}.}
\[\mathcal{C}^1_{\vartriangleright}=\bigcup_{r\in ]0,+\infty]}\mathcal{C}^1(]-\infty,r[),\]
une fonction $h\in \mathcal{C}^1_{\vartriangleright}$ représentant la réalisation $\mathfrak{h}$ de $\mathbb{H}$ définie sur 
\[D_{\mathfrak{h}}=D_h\cap \mathbf{R}_+^*\] 
par la formule (\ref{eq solution (D,f) de Histoire}) ci-dessus. Pour tout $h\in \mathcal{C}^1_{\vartriangleright}$, nous posons $\theta(h)=\mathfrak{h}$, où $\mathfrak{h}$ est la réalisation en question de $\mathbb{H}$. 

Pour constituer l'ensemble des réalisations de $\mathbb{H}$ pour lesquelles $D$ est fermé en sa borne supérieure, définissons  une copie disjointe de la partie 
\[C^{1\flat}_{\vartriangleright}=
\mathcal{C}^1_{\vartriangleright}
 \setminus \{h\in \mathcal{C}^1_{\vartriangleright}, D_h=]0,+\infty[\}
\]
de $\mathcal{C}^1_{\vartriangleright}$ formée des fonctions dont le domaine de définition est borné, copie que nous noterons $C^{1^\natural}_{\vartriangleright}$, et que nous définissons formellement en posant 
\[
C^{1^\natural}_{\vartriangleright}
= C^{1\flat}_{\vartriangleright}
\times \{\natural\},
\] où $\natural$ désignera également la bijection 
$\natural:C^{1\flat}_{\vartriangleright}
\ni h\mapsto 
h^\natural=(h,\natural)\in C^{1^\natural}_{\vartriangleright}$, ainsi que sa réciproque, de sorte que $(h,\natural)^\natural=h$ pour tout $h\in C^{1\flat}_{\vartriangleright}$. 

Nous décidons alors qu'un couple  $(h,\natural)\in C^{1^\natural}_{\vartriangleright}$ représentera la réalisation $\mathfrak{h}$ de $\mathbb{H}$ définie  par la formule (\ref{eq solution (D,f) de Histoire}), mais cette fois sur l'intervalle borné et fermé à droite \[D_\mathfrak{h}=\overline{D_h}\cap \mathbf{R}_+^*,\] et nous posons $\theta(h,\natural)=\mathfrak{h}$, où $\mathfrak{h}$ est la réalisation en question de $\mathbb{H}$.

Ainsi, l'ensemble des  réalisations\footnote{Des parties externes des réalisations, pour être précis (voir la remarque \ref{rmq usage realisation au lieu de partie externe} page \pageref{rmq usage realisation au lieu de partie externe}).} de $\mathbb{H}$ $=$ ${\mathbb{H}_{(k=1, a_0=-\infty, T_0=0, L={\mathcal{GC}_+^*})}}$ peut s'écrire comme l'union disjointe 
\[
\mathcal{S}_{\mathbb{H}}=
\{\underline{\emptyset}\}
\sqcup
\theta(\mathcal{C}^1_{\vartriangleright})
\sqcup
\theta(C^{1^\natural}_{\vartriangleright}),
\] 
avec $\underline{\emptyset}=\emptyset$. L'ensemble des réalisations non vides de $\mathbb{H}$ s'écrit quant à lui
\begin{equation}\label{eq ens real non vide de H}
\mathcal{S}_{\mathbb{H}}^*=
\theta(\mathcal{C}^1_{\vartriangleright})
\sqcup
\theta(C^{1^\natural}_{\vartriangleright}),
\end{equation}
et l'on a
\[\forall h\in \mathcal{C}^1_{\vartriangleright}, h=\widetilde{\theta(h)},\]
\[\forall h\in C^{1\flat}_{\vartriangleright}, h=\widetilde{\theta(h,\natural)}.\]

Pour toute réalisation non vide $\mathfrak{h}\in\mathcal{S}_{\mathbb{H}}^*$ de $\mathbb{H}$, nous appellerons \emph{partie mythique de l'histoire} la restriction $\widetilde{\mathfrak{h}}_{\vert ]-\infty,0]}$. Intuitivement, celle-ci représente en effet  une sorte d'histoire virtuelle, dont on aurait la mémoire sans qu'elle ait jamais été réellement vécue.

\begin{rmq}
Dans notre travail \cite{Dugowson:20170101} déjà cité,   la dynamique $\mathbb{H}_{(1, -\infty, 0, {\mathcal{GC}_+^*})}$ est notée \textcjheb{h} --- \emph{Hey}, en hébreu\footnote{En \TeX, nous codons la lettre \textcjheb{h} avec un \texttt{h} (code \texttt{\textbackslash textcjheb\{h\}}). }  ---  du fait de son rôle fondamental aux côtés des dynamiques \textcjheb{y} (exemple \ref{exm dyna historique}) et   \textcjheb{w} (exemple \ref{exm dynamique intemporelle}). 
\end{rmq}

\end{exm}

\begin{exm}[Autre exemple de dynamique de la forme $\mathbb{H}_{(k, a_0, T_0, L)}$]

Pour $k=0$,  $a_0=0$, $T_0=-\infty$ et $L=\{\mathbf{R}^2\}$, la dynamique  $\mathbb{K}=\mathbb{H}_{(k, a_0, T_0, L)}$ est de type
$[\dot{\pi}\overline{\delta}{\phi}\mathbf{R}_+]$, autrement dit c'est une mono-dynamique fonctorielle,  pluraliste et de moteur $\mathbf{R}_+$.
On vérifie aisément qu'une réalisation non vide $\mathfrak{h}$ de la dynamique $\mathbb{K}$ est 
\begin{itemize}
\item soit définie sur un intervalle de la forme $I=]-\infty,t[$ ou $I=]-\infty,t]$, avec $t\leq 0$, et dans ce cas $\mathfrak{h}(s)=(s,\emptyset)$ pour tout $s\in I$ (de sorte que $\mathfrak{h}$ est entièrement caractérisée par la donnée de $I$, à laquelle on peut l'identifier),
\item soit définie sur un intervalle de la forme $I=]-\infty,t[$ avec $t\in\mathbf{R}_+^*\cup\{+\infty\}$ ou de la forme $I=]-\infty,t]$ avec $t>0$, et on a $\mathfrak{h}(s)=(s,\emptyset)$ pour tout $s\in ]-\infty,0]$ et $\mathfrak{h}(s)=(s,h(s))$ pour tout $s>0$, avec $h\in\mathcal{C}(]0,t[)$ (et dans ce cas $\mathfrak{h}$ est entièrement caractérisée par la donnée de $I$ et de $h\in\mathcal{C}(\mathring{I}\cap \mathbf{R}_+^*)$).\\
\end{itemize}

Autrement dit, on a 
\[
\mathcal{S}_{\mathbb{K}}= \{\emptyset\}\cup
\left(
\bigcup_{t\in \mathbf{R}\cup\{+\infty\}}
(\{]-\infty,t[\}\times \mathcal{C}(]0,t[) )
\right)
\cup
\left(
\bigcup_{t\in \mathbf{R}}
(\{]-\infty,t]\}\times \mathcal{C}(]0,t[) )
\right).
\] 

\end{exm}


\subsubsection{Une dynamique intemporelle ($\mathbb{W}$ = \textcjheb{w})}\label{subsubs exm intemporel}

\begin{exm}[$\mathbb{W}$, une dynamique \og intemporelle \fg]\label{exm dynamique intemporelle}

Soit $L\subset\mathcal{C}$ un sous-ensemble non vide de  $\mathcal{C}$, 
et soit $\rightY$ une relation binaire entre  $L$ et $\mathcal{C}$.
Nous définissons la dynamique ${\mathbb{W}_\rightY}$, notée simplement $\mathbb{W}$ ci-après, en posant
\[
\mathbb{W}=
(
{\tau_{\mathbb{W}}}:
({\alpha_{\mathbb{W}}}:{\mathbf{C}_{\mathbb{W}}}\rightharpoondown \mathbf{P}^{\underrightarrow{\scriptstyle{L_{\mathbb{W}}}}})
\looparrowright 
({\mathbf{h}_{\mathbb{W}}}:{\mathbf{C}_{\mathbb{W}}}\rightarrow \mathbf{P})
),
\] 
avec
\begin{itemize}
\item pour moteur : ${\mathbf{C}_{\mathbb{W}}}=\mathbf{1}$, \emph{i.e.} $\dot{\mathbf{C}}_{\mathbb{W}}=\{\bullet\}$ et $\overrightarrow{\mathbf{C}_{\mathbb{W}}}=\{Id_\bullet=0\}$,
\item pour horloge, l'horloge existentielle $\mathbf{h}_{\mathbb{W}}=\xi_{\mathbf{1}}$,  qui n'a qu'un unique instant : $\bullet^{\mathbf{h}_{\mathbb{W}}}=\{0\}$ ,
\item pour paramétrisation $L_{\mathbb{W}}=
L$,
\item pour états : $st({\mathbb{W}})=\mathcal{C}$,
\item pour scansion la seule possible : $\forall f\in st(\mathbb{W}), \tau_{\mathbb{W}}(f)=0$,
\item pour tout $\omega\in L$, la transition $0^{\mathbb{W}}_\omega$ est définie pour toute $f\in st({\mathbb{W}})$ par
\[0^{\mathbb{W}}_\omega(f)=
\left[
\begin{array}{l}
\{f\}\mathrm{\quad si\,\,} \omega\rightY f, \\
\emptyset \mathrm{\quad sinon}.
\end{array} 
\right.
\]
\end{itemize}

Sauf dans le cas où la relation $\rightY$ est la relation complète, \footnote{Autrement dit de graphe  $L\times\mathcal{C}$, \emph{i.e.}  $\omega \rightY f$ est toujours vrai, et dans ce cas la dynamique $\mathbb{W}_\rightY$ est déterministe et donc fonctorielle.}, il y a toujours une valeur de $\omega$ pour laquelle au moins un état $f$ est \og hors-jeu \fg%
\footnote{Voir la section \textbf{§\,\ref{etats hors jeu}}.}, 
de sorte que la dynamique $\mathbb{W}_\rightY$ est de type  $[{\overline{\pi}}\underaccent{\dot}{\delta}\underline{\phi}\mathbf{1}]$, autrement dit paramétrique, bien quasi-déterministe, bien sous-fonctorielle et de moteur $\mathbf{1}=\{\bullet\}$.

Puisque l'horloge utilisée ne possède qu'un unique instant, (la partie externe d')une réalisation non vide s'identifie à un état de la dynamique. Précisons ce que sont ces réalisations pour deux choix particuliers de la relation $\rightY$.

\paragraph{Cas où $\rightY$ est la relation de compatibilité définie par :}

\[\forall (\omega,f)\in L\times \mathcal{C}, (\omega\rightY f)
\Leftrightarrow 
f_{\vert D_f\cap D_\omega}=\omega_{\vert D_f\cap D_\omega},\]
avec $L=\mathcal{C}$. Pour $\omega\in \mathcal{C}$, la mono-dynamique $\mathbb{W}_{\omega}=\mathbb{W}_{\rightY,\omega}$ a pour états \og dans le jeu\fg\, toutes les fonctions  $f\in\mathcal{C}$ définies et continues sur un intervalle ouvert $D_f$ qui coïncident avec $\omega$ sur $ D_f\cap D_\omega$, y compris la fonction vide.  Une réalisation non vide quelconque de $\mathbb{W}_{\omega}$ s'identifie naturellement à l'unique état par laquelle  elle passe en l'unique instant de l'horloge $\xi_{\mathbf{1}}$, donc à une fonction $f\in\mathcal{C}$, mais nous devons alors en particulier bien distinguer la réalisation non vide qui s'identifie à la fonction vide $\emptyset:\emptyset\hookrightarrow\mathbf{R}$, réalisation que nous noterons encore $\emptyset$, et la réalisation vide à proprement parler, $\underline{\emptyset}_\mathbb{W}$. Ainsi, concernant les réalisations de la dynamique $\mathbb{W}$, nous aurons
\[\underline{\emptyset}_\mathbb{W} \neq \emptyset\in \mathcal{S}_{\mathbb{W}}^*.\]

Par ailleurs, puisque $\emptyset\in\mathcal{F}$,  
nous pouvons en particulier considérer les réalisations non vides de la mono-dynamique ${\mathbb{W}_{\emptyset}}$, qui sont toutes les fonctions continues, de sorte que 
 $\mathcal{S}_{\mathbb{W}_{\emptyset}}=\{\underline{\emptyset}_\mathbb{W}\}\cup\mathcal{C}$, d'où, \emph{a fortiori},
\[\mathcal{S}_{\mathbb{W}}=\{\underline{\emptyset}_\mathbb{W}\}\cup\mathcal{C},\]
et pour l'ensemble des réalisations non vides de $\mathbb{W}$ :
\begin{equation}\label{eq ens real non vide de W}
\mathcal{S}_{\mathbb{W}}^*=
\mathcal{C}.
\end{equation}

\begin{rmq}
Cette dynamique est notée \textcjheb{w} --- \emph{vav}, en hébreu
\footnote{En \TeX, nous codons la lettre 
\textcjheb{w}
avec un
\texttt{w} 
(code \texttt{\textbackslash textcjheb\{w\}})
  . }
  --- dans notre travail \cite{Dugowson:20170101}, par référence au concept fondamental associé à la lettre \textcjheb{w} dans  la théorie \emph{métachrono\-logique} de Pierre Michel Klein \cite{KleinPM:2014}. 
\end{rmq}

\paragraph{Cas où $\rightY$ est l'injection canonique $L\hookrightarrow\mathcal{C}$ d'une partie $L$ de $\mathcal{C}$.}

Pour $\omega\in L$, la dynamique $\mathbb{W}_{\omega}=\mathbb{W}_{\hookrightarrow,\omega}$ a pour seul état \og dans le jeu\fg\, $f=\omega\in L\subset \mathcal{C}$. L'ensemble des (parties externes des) réalisations de la dynamique $\mathbb{W}=\mathbb{W}_{\hookrightarrow}$ s'identifie alors à $L$.

\end{exm}

\section{Familles interactives}\label{sec FI}

Le but principal de cette section est de définir et de donner un ou deux exemples de la notion de famille interactive, à savoir une famille de dynamiques ouvertes en interaction et coordonnées par des synchronisations. Après avoir précisé quelques notations, la notion d'interaction entre dynamiques ouvertes est précisée en section \textbf{§\,\ref{subs interactions dans une famille}}, celle de synchronisation en section \textbf{§\,\ref{subs synchronisation}}, la définition et des exemples de familles interactives étant proposés en section \textbf{§\,\ref{subs familles interactives}}.

\subsection{Notations}\label{subs Notations pour les familles interactives}

On se donne à partir de maintenant  une famille $\mathcal{A}=(A_i)_{i\in I}$ indexée par un ensemble non vide $I$ de dynamiques sous-fonctorielles ouvertes efficientes\footnote{Voir la définition \ref{df dynamique efficiente}.}
\[A_i=({\tau_i:(\alpha_i:\mathbf{C}_i\rightharpoondown \mathbf{P}^{\underrightarrow{\scriptstyle{L_i}}}})\looparrowright \mathbf{h}_i).\]
Pour tout $i\in I$, l'ensemble $\mathcal{S}^*_{A_i}$ des parties externes des réalisations non vides de la dynamique ouverte $A_i$ sera plus simplement noté
$\mathcal{S}^*_{i}=\mathcal{S}^*_{A_i}$, 
autrement dit
\[\mathcal{S}^*_i=
\bigcup_{\lambda\in L_i}
\mathcal{S}^*_{(i,\lambda)}, \]
 où $\mathcal{S}^*_{(i,\lambda)}=\mathcal{S}_{(i,\lambda)}\setminus \{ \underline{\emptyset}_{A_i}\}$ désigne l'ensemble des réalisations non vides de paramètre $\lambda$ de la dynamique $A_i$.
On notera en outre $\mathcal{Z}^*=(\mathcal{S}^*_{i})_{i\in I}$ la famille des $\mathcal{S}^*_i$, $\mathcal{L}=(L_i)_{i\in I}$ celle des $L_i$, et $\mathcal{E}=(E_i)_{i\in I}=(\mathcal{S}^*_i\times L_i)_{i\in I}$.

\subsection{Interactions dans une famille de dynamiques ouvertes efficientes}
\label{subs interactions dans une famille}

\subsubsection{Relations entre réalisations et paramètres}

\begin{df} 
Une \emph{relation entre réalisations et paramètres de la famille $\mathcal{A}=(A_i)_{i\in I}$} est une relation binaire multiple non vide  $C\in \mathcal{BM}^*_{(\mathcal{Z}^*, \mathcal{L})}$ de contexte d'entrée  
$\mathcal{Z}^*=(\mathcal{S}^*_{i})_{i\in I}$
et de contexte de sortie 
$\mathcal{L}=(L_i)_{i\in I}$.   
\end{df}

Le contexte global d'une telle relation entre réalisations et paramètres est alors $\mathcal{E}=(E_i)_{i\in I}=(\mathcal{S}^*_i\times L_i)_{i\in I}$, et nous écrirons souvent sous la forme 
$
\left(
\begin{array}{c}
\lambda_i  \\ 
\sigma_i 
\end{array} 
\right)_{i\in I}$
les éléments de $\Pi_I\mathcal{E}$  avec, pour tout $i\in I$,
 $\sigma_i\in \mathcal{S}^*_i$ et $\lambda_i \in L_i$. En particulier, un élément $\varsigma$ du graphe d'une relation $C$ entre réalisations et paramètres pour la famille $(A_i)_{i\in I}$ pourra s'écrire sous la forme
$\varsigma=
\left(
\begin{array}{c}
\lambda_i  \\ 
\sigma_i 
\end{array} 
\right)_{i\in I}
\in \vert C\vert
\subset
\Pi_I\mathcal{E}$.

\subsubsection{Cohérence, configurations et interactions}

\begin{df} Un
élément $\varsigma=
\left(
\begin{array}{c}
\lambda_i  \\ 
\sigma_i 
\end{array} 
\right)_{i\in I}\in\Pi_I\mathcal{E}$ est dit
\emph{cohérent pour la famille $\mathcal{A}$} s'il vérifie 
\[\sigma_i\in \mathcal{S}^*_{(i,\lambda_i)}\]
pour tout $i\in I$.
Un tel élément $\varsigma\in\Pi_I\mathcal{E}$ cohérent pour la famille $\mathcal{A}$ est également appelé une \emph{configuration pour  $\mathcal{A}$}.
\end{df}

Nous noterons provisoirement%
\footnote{Une nouvelle notation sera introduite en section \textbf{§\,\ref{subsubs interaction nulle}}.}
 ${CF}_\mathcal{A}$  l'ensemble des configurations pour $\mathcal{A}$.

\begin{df}\label{df interaction}[Interaction]
Une partie $G\subset\Pi_I\mathcal{E}$ est  dite \emph{cohérente} pour la famille $\mathcal{A}$ si chacun de ses éléments est cohérent, autrement dit si $G\subset {CF}_\mathcal{A}$.
Une relation $C\in \mathcal{BM}^*_{(\mathcal{Z}^*, \mathcal{L})}$ entre réalisations et paramètres pour la famille $\mathcal{A}$ est dite \emph{cohérente} si son graphe est cohérent, \textit{i.e.} si $\vert C\vert\subset {CF}_\mathcal{A}$, et une telle relation cohérente est également appelée une \emph{interaction dans $\mathcal{A}$}.
Étant donnée $R$ une interaction dans $\mathcal{A}$, on appelle \emph{$R$-configuration} toute configuration $\varsigma\in\vert R\vert$.
\end{df}

Nous noterons $\mathcal{T}_{\mathcal{A}}$ l'ensemble des interactions pour la famille $\mathcal{A}$. Toute interaction $R$ pour $\mathcal{A}$ étant caractérisée par la donnée d'un graphe non vide $\vert R\vert \subset {CF}_\mathcal{A}$, on en déduit une bijection canonique $\mathcal{T}_{\mathcal{A}}\simeq \mathcal{P}^*({CF}_\mathcal{A})$. 

\begin{rmq}
Soulignons le fait que, par définition, une dynamique inefficiente%
\footnote{Voir la définition \ref{df dynamique efficiente} page \pageref{df dynamique efficiente}.}
 ne peut entrer dans aucune interaction. Certes, comme c'était d'ailleurs le cas dans les définitions proposées dans \cite{Dugowson:20150807} et \cite{Dugowson:20150809},  nous aurions pu admettre des interactions vides ainsi que des interactions non vides mais contenant des réalisations vides, ce qui aurait nécessairement été le cas dès qu'une des dynamiques en jeu eût été inefficiente,  mais les dynamiques engendrées\footnote{Voir plus loin la section \ref{sec engendrement dynamique}.} par des familles interactives reposant sur de telles interactions auraient alors non seulement été elles-mêmes inefficientes, mais auraient en fait été totalement vides, \emph{i.e.} n'admettant pour seules transitions que les transitions vides. L'inefficience se révélant ainsi d'une stérilité absorbante, nous avons préféré la maintenir formellement en dehors de toute idée d'interaction.
\end{rmq}

\begin{rmq}\label{rmq compatibilite virtuelle des realisations}
Intuitivement, les configurations constituant une interaction expriment les comportements des diverses dynamiques qui sont mutuellement \og compatibles\fg\, lorsque ces dynamiques sont liés par cette interaction. Néan\-moins, une telle compatibilité peut s'avérer parfois quelque peu virtuelle, comme l'illustrera plus loin, page \pageref{rmq realisations non compatibles dans configurations},  la remarque \ref{rmq realisations non compatibles dans configurations} figurant dans l'exemple \ref{exm famille why} présenté en section \textbf{§\,\ref{subsubs famille why}}, remarque dans laquelle nous expliquons pourquoi cette virtualité éventuelle des compatibilités qu'exprime une interaction ne nous semble pas gênante.
\end{rmq}

Cela posé, il s'avère que la notion d'interaction donnée par la définition \ref{df interaction} est souvent trop générale, car elle permet des relations entre les dynamiques en jeu qui \og court-circuitent\fg\, le dispositif \emph{ad hoc} pour lequel l'influence mutuelle entre les dynamiques passe par les paramètres, comme l'illustre l'exemple \ref{exm interaction paranormale entre mono dynamiques} ci-dessous.

\begin{exm}\label{exm interaction paranormale entre mono dynamiques}
On considère la mono-dynamique fonctorielle ouverte de type $[\dot{\pi} \overline{\delta} \phi \mathbf{N}]$ :
\[A=({\tau:(\alpha:\mathbf{N} \rightarrow \mathbf{P}})\looparrowright \xi_\mathbf{N}),\]
définie par $st(\alpha)=\mathbf{N}\times\mathbf{R}$ et, pour tout $(n,r)\in st(\alpha)$,
\begin{itemize}
\item $\tau(n,r)=n$,
\item $\forall d\in\mathbf{N}^*, d^\alpha(n,r)=\{n+d\}\times \mathbf{R}$.
\end{itemize}

L'ensemble paramétrique de $A$ est un singleton, disons $L=\{*\}$. L'ensemble $\mathcal{S}^*$ de ses réalisations non vides s'identifie à l'ensemble des suites finies ou infinies de réels $\sigma=(s_n)_{n\in E_\sigma\subset \mathbf{N}}$, avec $E_\sigma$ un intervalle de $\mathbf{N}$ qui contient $0$. 
Considérons alors le cas où la famille $\mathcal{A}=(A_i)_{i\in I}$ est définie par $I=\{1,2\}$ et $A_1=A_2=A$, et soit $R$ la relation entre réalisations non vides et paramètres pour cette famille définie par 
\[ \forall (\sigma_1,\sigma_2)\in\mathcal{S}^*_1\times\mathcal{S}^*_2=(\mathcal{S}^*)^2,
\left(
\begin{array}{cc}
 * & * \\ 
 \sigma_1 & \sigma_2
\end{array}  
\right)
\in \vert R\vert 
\Leftrightarrow
\sigma_1=\sigma_2.
\]
Cette relation est trivialement cohérente, c'est donc bien une interaction. Cependant, cette interaction consiste en une relation directe entre les réalisations des dynamiques en jeu, et ne fait jouer aucun rôle aux paramètres, et pour cause : nous avons ici affaire à des mono-dynamiques, qui ne sont pas censées être influençables. Une telle interaction ne sera pas considérée comme \og normale\fg. Ceci sera précisé en section \textbf{§\,\ref{subsubs interactions normales}}, grâce en particulier aux notions de \emph{partie cohérente d'une relation} (section \textbf{§\,\ref{subsubs partie coherente de RRP}}) et de \emph{relation filtrante} (section \textbf{§\,\ref{subsubs relations filtrantes et interactions operantes}}).
\end{exm}

\subsubsection{Partie cohérente d'une relation entre réalisation et paramètres}\label{subsubs partie coherente de RRP}

Nous définissons une application 
\[\mathcal{BM}^*_{(\mathcal{Z}^*, \mathcal{L})}\ni C
\mapsto 
\widecheck{C} \in
\mathcal{T}_{\mathcal{A}}\cup\{\emptyset\}
\subset\mathcal{BM}_{(\mathcal{Z}^*, \mathcal{L})}
\] 
en associant à toute relation  entre réalisations et paramètres $C\in \mathcal{BM}^*_{(\mathcal{Z}^*, \mathcal{L})}$ de la famille $\mathcal{A}$ la relation multiple $\widecheck{C}\in \mathcal{BM}_{(\mathcal{Z}^*, \mathcal{L})}$ définie par son graphe
\[
\vert \widecheck{C}\vert=
\vert C\vert\cap {CF}_\mathcal{A}.
\]

La relation multiple $\widecheck{C}$ 
sera appelée la \emph{partie cohérente pour $\mathcal{A}$ de ${C}$}. Si  $\widecheck{C}$ est non vide, c'est une interaction dans $\mathcal{A}$.

\subsubsection{L'interaction nulle et l'ensemble des configurations de $\mathcal{A}$}\label{subsubs interaction nulle}

Soit $TL$ la relation binaire multiple totale  de contexte $(\mathcal{Z}^*,\mathcal{L})$, autrement dit  la relation binaire multiple de graphe  $\vert TL\vert=\Pi_I\mathcal{E}$. 

\begin{df} On appelle \emph{interaction nulle pour $\mathcal{A}$}, et l'on note $\Omega_\mathcal{A}$, la partie cohérente de $TL$ :
\[\Omega_\mathcal{A}=\widecheck{TL}.\]
\end{df}

Le graphe de l'interaction nulle $\Omega_{\mathcal{A}}$ étant par définition 
$
\vert\Omega_{\mathcal{A}}\vert={CF}_\mathcal{A},
$  
nous pourrons à l'avenir noter $\vert\Omega_{\mathcal{A}}\vert$ plutôt que ${CF}_\mathcal{A}$ l'ensemble des configurations de $\mathcal{A}$, et c'est ce que nous ferons.

\begin{rmq} L'interaction nulle est  l'interaction \emph{maximale} du point de vue de l'inclusion des graphes. Inversement, plus le graphe d'une interaction est mince, plus on pourra considérer que l'interaction est forte entre les dynamiques en jeu : interagir, c'est restreindre les possibles\footnote{C'est d'ailleurs vrai également pour toute création, ce qu'illustre admirablement la notion de contrainte créatrice à l'\oe uvre dans la littérature oulipienne.}. \emph{A contrario}, quand \og tout est possible\fg, c'est qu'il n'y a en réalité pas vraiment d'interaction.
\end{rmq}

\subsubsection{Relations filtrantes et interactions opérantes}\label{subsubs relations filtrantes et interactions operantes}

\paragraph{Ensemble $\mathcal{X}(R)$.}

Pour toute interaction $R\in\mathcal{T}_\mathcal{A}$, nous posons
\[\mathcal{X}(R)=\{C\in\mathcal{BM}^*_{(\mathcal{Z}^*, \mathcal{L})}, \widecheck{C}=R\}.\]

\begin{rmq}\label{rmq interpretation ensemble XR}
Les relations entre réalisations et paramètres $C\in\mathcal{X}(R)$ ne sont pas, en général, cohérentes, mais peuvent être interprétées intuitivement comme exprimant une \og demande\fg\, d'interaction, la \og réponse\fg\, cohérente à cette demande étant précisément l'interaction $R$. Dès lors, nous pouvons voir intuitivement l'ensemble $\mathcal{X}(R)$ comme celui des interprétations possibles de l'interaction $R$ en termes de demandes non nécessairement cohérentes. 
\end{rmq}

\begin{df} Une relation entre réalisations et paramètres $C\in \mathcal{BM}^*_{(\mathcal{Z}^*, \mathcal{L})}$ pour la famille $\mathcal{A}$ est dite \emph{filtrante} si le domaine de définition de la relation binaire qu'elle détermine $\Pi_I\mathcal{Z}^*\rightarrow \Pi_I\mathcal{L}$ n'est pas $\Pi_I\mathcal{Z}^*$ tout entier, autrement dit si
\[
D_{rb(C)}\subsetneqq \Pi_I\mathcal{Z}^*.
\]
\end{df}

\begin{df}
Une interaction $R\in\mathcal{T}_{\mathcal{A}}$ est dite \emph{opérante} si elle est filtrante.  \emph{A contrario}, une interaction non filtrante sera dite \emph{inopérante}. 
\end{df}

\begin{exm}
L'interaction nulle $\Omega_{\mathcal{A}}$ est inopérante.
\end{exm}

\subsubsection{Structure connective des réalisations d'une interaction}

Intuitivement, une interaction inopérante exprime une absence d'interaction effective, puisque dans ce cas tous les comportements (les réalisations) des dynamiques en jeu sont compatibles. Ceci peut être précisé grâce à la notion de \emph{structure connective des réalisations d'une interaction}, dont la définition \ref{df structures connectives de interaction} ci-dessous s'appuie sur celle de \emph{structure connective d'une relation multiple} que nous avons introduite dans \cite{Dugowson:201505}. La définition \ref{df structures connectives de interaction} ci-dessous reprend également celle de structure connective d'une \og famille dynamique\fg\, introduite dans \cite{Dugowson:20150809}, mais en la désignant désormais comme la structure connective \emph{globale} de l'interaction concernée, pour la distinguer de la structure connective \emph{des réalisations} de cette même interaction.

\begin{df}\label{df structures connectives de interaction} Étant donnée $R\in\mathcal{T}_{\mathcal{A}}$  une interaction dans la famille de dynamiques sous-fonctorielles efficientes $\mathcal{A}=(A_i)_{i\in I}$ d'ensembles paramétriques $L_i$ et d'ensembles de réalisations non vides $\mathcal{S}^*_i$,
\begin{itemize}
\item la \emph{structure connective globale de l'interaction $R$} est la structure connective (sur $I$)\footnote{Concernant la notion de structure connective d'une relation multiple, voir la définition 17, section \textbf{§\,3.1} de \cite{Dugowson:201505}.} de la relation multiple\footnote{Concernant la relation multiple notée $rm(R)$, voir plus haut la section \textbf{§\,\ref{subsubs relations binaires multiples}}.}
 $rm(R)$ de contexte $\mathcal{E}={(\mathcal{S}^*_i\times L_i)}_{i\in I}$,
\item la \emph{structure connective des réalisations de l'interaction $R$} est la structure connective (sur $I$) de la relation multiple $\underline{R}$ de contexte $\mathcal{Z}^*=(\mathcal{S}^*_i)_{i\in I}$ définie par 
\[
(\sigma_i)_{i\in I}\in\underline{R} \Leftrightarrow \exists (\lambda_i)_{i\in I}\in \prod_{i\in I}L_i, 
\left(
\begin{array}{c}
\lambda_i  \\ 
\sigma_i 
\end{array} 
\right)_{i\in I}\in \vert R\vert.
\]
\end{itemize}
\end{df}

\begin{rmq}
La structure connective des réalisations d'une interaction est toujours plus fine (ou égale) que sa structure connective globale.
\end{rmq}

\begin{prop}
La structure connective sur $I$ des réalisations d'une interaction inopérante est discrète\footnote{Une structure connective discrète est également dite \emph{totalement déconnectée} : les seules parties connexes non vides de $I$ sont les singletons; voir \cite{Dugowson:201012}.}.
\end{prop}
\paragraph{Preuve.}
La relation multiple $\underline{R}$, projection de $R$ sur $\Pi_I \mathcal{Z}^*=\prod_{i\in I}{\mathcal{S}^*_i}$, étant la relation multiple totale sur $\Pi_I \mathcal{Z}^*$, la proposition ci-dessus résulte immédiatement de  la notion de structure connective d'une relation multiple telle qu'elle est définie dans \cite{Dugowson:201505}, section \textbf{§\,3.1}.
\begin{flushright}$\square$\end{flushright} 
\pagebreak[2]

\subsubsection{Interactions normales}\label{subsubs interactions normales}

\begin{df} Une interaction $R\in\mathcal{T}_{\mathcal{A}}\subset\mathcal{BM}^*_{(\mathcal{Z}^*, \mathcal{L})}$ pour la famille de dynamiques $\mathcal{A}$ est dite
\begin{itemize}
\item \emph{normale}, s'il existe $C\in\mathcal{X}(R)$ telle que $C$ ne soit pas filtrante,
\item \emph{déterminante}, si $rb(R)$ est une fonction\footnote{Voir les rappels de la section \textbf{§\,\ref{subsubs relations binaires applications fonctions}}.},
\item \emph{concrète} si elle est à la fois normale et déterminante.
\end{itemize}
\end{df}

\begin{rmq}\label{rmq interpretation interaction operante} Une interaction inopérante est normale. D'un autre côté, 
 si elle est opérante, une interaction normale $R$ filtre certes les réalisations des dynamiques en jeu, non pas de façon arbitraire et \emph{a priori}, mais du seul fait de l'exigence de cohérence. En effet, selon la remarque \ref{rmq interpretation ensemble XR} concernant l'interprétation intuitive de l'ensemble $\mathcal{X}(R)$, il existe pour une telle interaction $R$ une \og interprétation en termes de demande\fg\, non filtrante $C\in\mathcal{X}(R)$  de $R$, une telle demande n'ayant pas vocation à être cohérente. L'exemple \ref{exm interaction paranormale entre mono dynamiques} est celui d'une interaction qui n'est pas normale. Une telle interaction sera dite \emph{paranormale}, ce terme suggérant que les corrélations entre les réalisations des dynamiques en jeu ne passent pas nécessairement par les dispositifs paramétriques qui représentent en quelque sorte les \og organes sensoriels\fg\, de ces dynamiques.
\end{rmq}

\begin{rmq}
S'il existe $C\in\mathcal{X}(R)$ telle que $rb(C)$ soit une application de $\Pi_I\mathcal{Z}^*$ vers $\Pi_I\mathcal{L}$, alors $R$ est concrète. On vérifie facilement que, \emph{modulo} l'axiome du choix, l'existence d'une telle relation $C$ caractérise les interactions concrètes.
\end{rmq}

\subsection{Synchronisations}\label{subs synchronisation}

\begin{df}\label{df synchronisation}
Étant données $A_0$ et $A_1$ deux dynamiques sous-fonctorielles ouvertes, avec pour $i\in\{0,1\}$,
\[A_i=({\tau_i:(\alpha_i:\mathbf{C}_i\rightharpoondown \mathbf{P}^{\underrightarrow{\scriptstyle{L_i}}}})\looparrowright \mathbf{h}_i),\]
on appelle \emph{synchronisation de $A_1$ par $A_0$} tout couple $(\Delta_1,\delta_1)$ constitué
\begin{itemize}
\item d'une application $\Delta_1:\dot{\mathbf{C}}_0\rightarrow \dot{\mathbf{C}}_1$ associant aux \emph{objets} de $\mathbf{C}_0$ des objets de $\mathbf{C}_1$,
\item d'une application  $\delta_1:st(\mathbf{h}_0)\rightarrow st(\mathbf{h}_1)$ associant à tout instant de l'horloge $\mathbf{h}_0$ un instant de l'horloge $\mathbf{h}_1$ et qui vérifie les deux conditions suivantes : 
\subitem $-$
$\forall S\in\dot{\mathbf{C}}_0, \forall s\in S^{\mathbf{h}_0}, \delta_1(s)\in (\Delta_1 S)^{\mathbf{h}_1}$,
\subitem $-$ $\delta_1$ est monotone%
\footnote{Autrement dit $\delta_1$ est soit croissante, ce qui signifie que $s_0\leq_{\mathbf{h}_0} t_0\Rightarrow \delta_1(s_0)\leq_{\mathbf{h}_1} \delta_1(t_0)$, où $\leq_{\mathbf{h}_i}$ désigne le pré-ordre sur les instants de l'horloge $\mathbf{h}_i$ (voir la section \ref{subs instants et anteriorite}), soit décroissante : $s_0\leq_{\mathbf{h}_0} t_0\Rightarrow \delta_1(t_0)\leq_{\mathbf{h}_1} \delta_1(s_0)$.}.
\end{itemize}
Dans le cas particulier où $(\Delta_1,\delta_1)$ est un dynamorphisme
\footnote{
Plus exactement, $\Delta_1$ est alors la partie objet de la partie fonctorielle d'un tel dynamorphisme (voir plus haut page \pageref{paragraph categorie des horloges} le paragraphe \emph{catégorie des horloges} de la section \textbf{§\,\ref{paragraph categorie des horloges}}).}
, alors nécessairement déterministe, de $\mathbf{h}_0$ vers $\mathbf{h}_1$, nous dirons que la synchronisation $(\Delta_1,\delta_1)$ est \emph{rigide}. Bien entendu, dans le cas contraire, nous dirons que cette synchronisation est \emph{souple}.
\end{df}

Nous écrirons $(\Delta_1,\delta_1):\mathbf{h}_0\Rsh \mathbf{h}_1$ pour indiquer que $(\Delta_1,\delta_1)$ est une synchronisation de $\mathbf{h}_1$ par $\mathbf{h}_0$.
 
\begin{rmq}
La notion de synchronisation donnée par la définition \ref{df synchronisation} ci-dessus est considérablement plus générale que celle définie dans \cite{Dugowson:20150807} et \cite{Dugowson:20150809}, où seules des synchronisations rigides avaient été considérées.
\end{rmq}

\subsection{Familles interactives}\label{subs familles interactives}

\subsubsection{Définition}

\begin{df}\label{df famille interactive} On appelle \emph{famille interactive}  la donnée $(I,\mathcal{A}, R,i_0,  (\Delta_i,\delta_i)_{i\neq i_0})$
\begin{itemize}
\item d'un ensemble $I$ non vide, appelé \emph{index} de la famille,
\item d'une famille indexée par $I$ de dynamiques sous-fonctorielles ouvertes efficientes $\mathcal{A}=(A_i)_{i\in I}$, avec 
\[A_i=({\tau_i:(\alpha_i:\mathbf{C}_i\rightharpoondown \mathbf{P}^{\underrightarrow{\scriptstyle{L_i}}}})\looparrowright \mathbf{h}_i),\]
appelées \emph{composantes} de la famille interactive, 
\item d'une interaction $R\in \mathcal{T}_{\mathcal{A}}$ pour la famille $\mathcal{A}$, 
\item d'un élément $i_0\in I$, appelé \emph{indice synchronisateur} de la famille,
\item d'une famille 
\[
((\Delta_i,\delta_i) 
: \mathbf{h}_{i_0}
\Rsh
\mathbf{h}_i)_{i\in I\setminus\{i_0\}}
\] 
de synchronisations des $\mathbf{h}_i$ par $\mathbf{h}_{i_0}$.
\end{itemize}
\end{df}

\begin{rmq}
L'expression \og familles dynamiques\fg\, utilisée dans \cite{Dugowson:20150807} et \cite{Dugowson:20150809} pour les désigner présentait l'inconvénient d'une trop grande ressemblance avec celle de \og familles \emph{de} dynamiques\fg, alors que celles-ci ne constituent qu'une part de celles-là. Pour cette raison, nous préférons utiliser désormais l'expression \emph{familles interactives}.
\end{rmq}

\begin{rmq}
Par \emph{chef d'orchestre} d'une famille interactive, nous entendrons selon les contextes soit l'indice synchronisateur de cette famille, soit la dynamique ouverte dont l'indice dans la famille est l'indice synchronisateur.
\end{rmq}

\subsubsection{Structures connectives d'une famille interactive}

Remarquons au passage --- les questions connectives devant faire l'objet de publications ultérieures --- que de la définition \ref{df structures connectives de interaction} des deux structures connectives (respectivement \og globale\fg\, et \og des réalisations\fg) associées à une interaction découle immédiatement celles d'une famille interactive :

\begin{df} La \emph{structure connective globale (resp. des réalisations) d'une famille interactive} est la structure connective globale (resp. des réalisations) de son interaction\footnote{Voir plus haut la définition \ref{df structures connectives de interaction}.}. L'\emph{ordre
connectif global (resp. des réalisations) d'une famille interactive} est l'ordre connectif\footnote{Sur la notion d'ordre connectif, voir par exemple, dans le cas fini, la définition 16 de \cite{Dugowson:201012}, et dans le cas général, la section \textbf{§\,1.11} de \cite{Dugowson:201112} ou de  \cite{Dugowson:201203}.} de la structure connective globale (resp. des réalisations) de cette famille.
\end{df}

\subsubsection{Exemple de la famille $\mathbb{WHY}$ = \textcjheb{why} }\label{subsubs famille why} 
\begin{exm}\label{exm famille why} Poursuivons  l'exemple qui sera présenté dans \cite{Dugowson:20170101} à l'occasion de notre travail en collaboration avec le philosophe Pierre-Michel Klein sur la base de sa théorie \emph{métachronologique} \cite{KleinPM:2014}. Nous y proposons de considérer une famille interactive notée $\mathbb{WHY}$ ou  \textcjheb{why}, dont les composantes sont les trois dynamiques \textcjheb{y} $ = \mathbb{Y}$ (source), \textcjheb{h} $ = \mathbb{H}$ (histoire), et \textcjheb{w} $ =\mathbb{W}$  (intemporel)  présentées précédemment (exemples \ref{exm source lip}, \ref{exm dyna historique} et \ref{exm dynamique intemporelle}). Précisément, reprenant les notations introduites dans ces exemples\footnote{En particulier les applications $\theta$ et $\mathfrak{h}\mapsto \widetilde{\mathfrak{h}}$ introduites dans l'exemple \ref{exm dyna historique}.} nous posons
\[
\mathbb{WHY}=\textcjheb{why}=
(I=\{1,2,3\}, \mathcal{A}=(A_1, A_2, A_3), R,i_0=2, (\Delta_i,\delta_i)_{i\in\{1,3\}}),
\] 
avec
\begin{itemize}
\item $A_1 =\mathbb{Y}= $ \textcjheb{y}, la dynamique \og source\fg\, donnée à l'exemple \ref{exm source lip},
\item $A_2 =\mathbb{H}=\mathbb{H}_{(k=1, a_0=-\infty, T_0=0, L={\mathcal{GC}_+^*})}= $   \textcjheb{h},  la dynamique \og historique\fg\, de l'exemple \ref{exm dyna historique},
\item $A_3 =\mathbb{W} = $ \textcjheb{w},  la dynamique \og intemporelle\fg\, de l'exemple \ref{exm dynamique intemporelle}, dans le cas où la relation de compatibilité est $(\lambda\rightY f)
\Leftrightarrow 
f_{\vert D_f\cap D_\lambda}=\lambda_{\vert D_f\cap D_\lambda}$,
\item $i_0=2$, autrement dit le chef d'orchestre est $\mathbb{H}= $  \textcjheb{h},
\item $\Delta_1=Id_{\mathbf{R}_+}$ 
et $\delta_1$ est l'injection canonique 
$st(\mathbf{h}_\textcjheb{h})=]0,+\infty[$ 
$\hookrightarrow$  
$[0,+\infty[=st(\mathbf{h}_\textcjheb{y})$,
\item  $\Delta_3=(\mathbf{R}_+ \rightarrow \mathbf{0})$ et $\delta_3:st(\mathbf{h}_\textcjheb{h})=]0,+\infty[\rightarrow \{0\}=st(\mathbf{h}_\textcjheb{w})$ sont, par nécessité évidente, constants,
\item enfin, on prend pour interaction $R\in \mathcal{T}_{\mathcal{A}}$ celle dont le graphe $\vert R\vert$ est l'ensemble des configurations
\[\left(
\begin{array}{ccc}
\omega \in \mathcal{C}		&
\gamma \in \mathcal{GC}_+^*		& 
* \\ 
\mathfrak{w}\in\mathcal{C}	& 
\mathfrak{h}\in \theta(\mathcal{C}^1_{\vartriangleright}
\sqcup
C^{1^\natural}_{\vartriangleright}) 	& 
\mathfrak{y}\in (Lip^1_+)^*
\end{array}	 
\right)\in \vert \Omega_{\mathcal{A}}\vert\]
qui vérifient
$\gamma=\mathfrak{y}_{\vert int(D_\mathfrak{y}) }$
et
$\omega=\widetilde{\mathfrak{h}}$.\\
\end{itemize}

Ainsi, $\vert R\vert$ est  l'ensemble des 
$\left(
\begin{array}{ccc}
\omega 	&
\gamma 		& 
* 			\\ 
\mathfrak{w} & 
\mathfrak{h} & 
\mathfrak{y} 
\end{array}	 
\right)\in \Pi_\mathcal{E}$ tels que
\begin{itemize}
\item $\mathfrak{y}\in (Lip^1_+)^*$,
\item $\gamma=\mathfrak{y}_{\vert int(D_\mathfrak{y})}\in\mathcal{GC}_+^*$,
\item $\widetilde{\mathfrak{h}}$ est définie et de classe $\mathcal{C}^1$ sur un intervalle de la forme $]-\infty,r[$ et coïncide avec $\gamma$ sur $]0,r[= D_\gamma$,
\item $\omega=\widetilde{\mathfrak{h}}\in \mathcal{C}$.
\item $\omega \rightY  \mathfrak{w}\in\mathcal{C}$.
\end{itemize}

Puisque $\mathfrak{y}$ détermine $\gamma$ et que $\mathfrak{h}$ détermine $\omega$, il est clair que l'interaction $R$ est concrète.

\begin{rmq} Des relations entre les différentes composantes d'une $R$-configuration quelconque
$\left(
\begin{array}{ccc}
\omega 	&
\gamma 		& 
* 			\\ 
\mathfrak{w} & 
\mathfrak{h} & 
\mathfrak{y} 
\end{array}	 
\right)\in \vert R\vert$, on déduit les propriétés suivantes :
\begin{itemize}
\item il existe nécessairement $s>0$ tel que $\mathfrak{y}_{\vert [0,s[}\in\mathcal{C}^1([0,s[)$, et l'ensemble des $s$ vérifiant cette propriété est un intervalle de borne supérieure $s_{\max}\in \overline{\mathbf{R}}_+^*$,
\item  $\gamma_{\vert ]0,s_{\max}[}\in\mathcal{C}^1(]0,s_{\max}[)$ et $\gamma\in (Lip^1_+)^*$,
\item $\exists t\in ]0,s_{\max}]$, ($D_\mathfrak{h}=]0,t[$ ou $D_\mathfrak{h}=]0,t]$) et $\widetilde{\mathfrak{h}}_{\vert ]0,t[}=\mathfrak{y}_{\vert ]0,t[}$,
\item $\omega=\widetilde{\mathfrak{h}}$ donc $\omega\in \mathcal{C}^1_{\vartriangleright}$,
\item $\mathfrak{w}_{\vert D_\mathfrak{w} \cap ]-\infty,t[}$ est de classe $\mathcal{C}^1$ par compatibilité avec $\omega$.
\end{itemize}
\end{rmq}

\begin{rmq}\label{rmq realisations non compatibles dans configurations}
Illustrant la remarque \ref{rmq compatibilite virtuelle des realisations} faite plus haut (page \pageref{rmq compatibilite virtuelle des realisations}), notons que dans une $R$-configuration $\left(
\begin{array}{ccc}
\omega 	&
\gamma 		& 
* 			\\ 
\mathfrak{w} & 
\mathfrak{h} & 
\mathfrak{y} 
\end{array}	 
\right)$, $\mathfrak{w}$ et $\mathfrak{y}$ \emph{ne sont pas nécessairement compatibles} (au sens de la relation $\rightY$). Par exemple, si $f:]-\infty,4]\rightarrow \mathbf{R}$ est une fonction qui est de classe $\mathcal{C}^1$ sur $]-\infty,3[$, qui est $1$-lipschitzienne sur $[0,4]$ et qui est strictement croissante sur $[2,3]$, alors en posant $\mathfrak{y}=f_{\vert [0,4]}$, $\gamma=f_{\vert ]0,4[}$, $\mathfrak{h}=\omega=f_{\vert ]-\infty,2[}$ et $\mathfrak{w}:[1,3]\rightarrow \mathbf{R}$ tel que $\mathfrak{w}(x)= f(x)$ si $x\leq 2$ et $\mathfrak{w}(x)= f(2)$ si $x\geq 2$, on obtient une $R$-configuration $\left(
\begin{array}{ccc}
\omega 	&
\gamma 		& 
* 			\\ 
\mathfrak{w} & 
\mathfrak{h} & 
\mathfrak{y} 
\end{array}	 
\right)$ telle que pour tout $x\in]2,3]$ on a $\mathfrak{w}(x)\neq\mathfrak{y}(x)$. 

On pourrait considérer que l'existence de telles réalisations non compatibles est un défaut des interactions telles que nous les avons définies, mais ce n'est pas notre sentiment, d'une part car ce \og défaut\fg\, ne gênera en rien, ci-après dans la section \textbf{§\,\ref{sec engendrement dynamique}}, la définition des dynamiques engendrées par une famille interactive, d'autre part et surtout parce qu'une définition générale de la compatibilité des réalisations ne pourrait que s'appuyer sur ces dynamiques engendrées et apporterait de ce fait des complications importantes et inutiles à la notion d'interaction. Il nous paraît donc préférable d'admettre qu'au sein des interactions des configurations en quelques sorte virtuelles, comportant des réalisations incompatibles, puisse apparaître, quitte à réfléchir aux conséquences philosophiques de l'existence de telles configurations.
\end{rmq}

\begin{rmq}
Nous verrons ultérieurement plusieurs dynamiques produites par la famille interactive $\mathbb{WHY}$ = \textcjheb{why}, en particulier, en  section \textbf{§\,\ref{subsubs exemple shin}}, celle que nous noterons $\mathbb{S}$ ou \textcjheb{/s}.
\end{rmq}

\end{exm}

\section{Engendrement dynamique}\label{sec engendrement dynamique}

\begin{rmq}\label{rmq ordre des facteurs}
Étant donné $I$ l'ensemble non vide  indexant une famille interactive $\mathcal{F}$, nous  considérerons que
l'ordre des facteurs n'intervient pas dans les produits cartésiens d'ensembles indexés par $I$ (ou, dans certains cas, par $2I$) que nous écrirons dans cette section, chaque facteur étant en tout état de cause associé sans ambiguïté à un indice précis $i\in I$ (ou, le cas échéant, à $j\in 2I$). Ainsi, pour une famille d'ensemble $(T_i)_{i\in I}$ et un élément $i_0\in I$, les expressions $T_{i_0}\times \prod_{i\neq {i_0}}{T_i}$ et $(\prod_{i\neq {i_0}}{T_i})\times T_{i_0} $ devront être comprises comme désignant toutes deux $\prod_{i\in I}{T_i}$ :
\[
T_{i_0}\times \prod_{i\neq {i_0}}{T_i}
=(\prod_{i\neq {i_0}}{T_i})\times T_{i_0} 
=\prod_{i\in I}{T_i}.
\]
\end{rmq}

\subsection{Théorème de stabilité sous-fonctorielle}\label{subs thm stab}

\begin{prop}\label{prop multidynamique graphique est sous-fonctorielle} Soit $\mathcal{F}=(I,\mathcal{A}=(A_i)_{i\in I}, R, i_0, (\Delta_i,\delta_i)_{i\neq i_0})$ une famille interactive, de composantes les dynamiques sous-fonctorielles ouvertes efficientes
\[A_i=({\tau_i:(\alpha_i:\mathbf{C}_i\rightharpoondown \mathbf{P}^{\underrightarrow{\scriptstyle{L_i}}}})\looparrowright \mathbf{h}_i),\]
et soit $\beta_\mathcal{F}=\beta:\vert\mathbf{B}\vert\longrightarrow \vert\mathbf{P}^{\underrightarrow{\scriptstyle{M}}}\vert$ 
la multi-dynamique graphique 
définie par  \\

\begin{itemize}
\item $\mathbf{B}=\mathbf{C}_{i_0}$,
\item $M=Im(rb(R))$,
\item pour tout sommet $S\in \dot{\vert\mathbf{B}\vert} =\dot{\mathbf{B}}$,
\[S^\beta
=
\{(a_i)_{i\in I}\in S^{\alpha_{i_0}}\times \prod_{i\neq {i_0}}(\Delta_i S)^{\alpha_i},
\forall i\neq {i_0}, \tau_{i}(a_i)=\delta_i(\tau_{{i_0}}(a_{i_0}))\},\]
\item pour toute arête $(e:S\rightarrow T)\in \overrightarrow{\vert\mathbf{B}\vert}=\overrightarrow{\mathbf{B}}$, tout état  $a=(a_i)_{i\in I}\in S^\beta$ et tout paramètre $\mu\in M$,
\[ 
e^\beta_\mu(a)=
\{
b=(b_i)_{i\in I}\in T^\beta, \,
b \mathrm{\,v\acute{e}rifie\,les\,conditions\,\,} (\ref{eq 1 sync})\mathrm{\,et\,}(\ref{eq 2 real})
\} 
\]les conditions en question étant respectivement
\begin{equation}\label{eq 1 sync}
\tau_{{i_0}}(b_{i_0})=e^{\mathbf{h}_{i_0}}(\tau_{{i_0}}(a_{i_0}))
\end{equation}
et
\begin{equation}\label{eq 2 real}
\exists (\mathfrak{a}_i)_{i\in I}\in  rb(R)^{-1}(\mu),
\forall i\in I, \mathfrak{a}_i\triangleright a_i, b_i.
\end{equation}\\
\end{itemize}
Alors $\beta$ est une multi-dynamique sous-fonctorielle.
\end{prop}
\paragraph{Preuve.} Remarquons tout d'abord que $\beta$ est bien une multi-dynamique graphique telle que définie en section \ref{subsubs multi dynamiques graphiques}, puisque pour $(S,T)\in\dot{\mathbf{B}}^2$ avec $S\neq T$ on a 
$\forall (a_i)_{i\in I}\in S^\beta\cap T^\beta,  a_{i_0}\in S^{\alpha_{i_0}}\cap T^{\alpha_{i_0}}=\emptyset$ de sorte que $S^\beta\cap T^\beta=\emptyset$.

Ensuite, conformément à la remarque \ref{rmq condition graphique soit sous fonct} (page \pageref{rmq condition graphique soit sous fonct}), nous devons vérifier les deux conditions suivantes :
 
\begin{itemize}
\item 
$
\forall S\in\dot{\mathbf{B}}, \forall \mu\in M, (Id_S)^\beta_\mu\subset Id_{S^\beta}
$,
\item 
$
\forall (S{\stackrel{f}{\rightarrow}}T{\stackrel{g}{\rightarrow}}U)\in\overrightarrow{\mathbf{B}}^2,
\forall \mu\in M, 
(g\circ f)^\beta_\mu\subset g^\beta_\mu \odot f^\beta_\mu
$.
\end{itemize}

\subparagraph{Vérifions la première condition.} Soit donc  $a={(a_i)_{i\in I}}\in S^\beta$. Supposons $(Id_S)^\beta_\mu(a)\neq \emptyset$, et soit  $a'={(a'_i)_{i\in I}}\in (Id_S)^\beta_\mu(a)$. Pour tout $i\in I$, il existe $\mathfrak{a}_i\in\mathcal{S}^*_{i}=\mathcal{S}^*_{A_i}$ telle que $\mathfrak{a}_i\rhd a_i,a'_i$  et l'on a, par définition de $S^\beta$ et par la condition (\ref{eq 1 sync}) appliquée à $e=Id_S$,
\[
a'_i
= \mathfrak{a}_i(\tau_i(a'_i))
=\mathfrak{a}_i(\delta_i(\tau_{i_0}(a'_{i_0})))
=\mathfrak{a}_i(\delta_i(\tau_{i_0}(a_{i_0})))
= \mathfrak{a}_i(\tau_i(a_i))
=a_i,
\]
de sorte que dans tous les cas $(Id_S)^\beta_\mu(a)\subset \{a\}$, autrement dit, en vertu de la remarque \ref{rmq appli det} (page \pageref{rmq appli det}), $(Id_S)^\beta_\mu(a)\subset Id_{S^\beta}(a)$, d'où
\[(Id_S)^\beta_\mu\subset Id_{S^\beta}.\]

\subparagraph{Vérifions à présent la seconde condition.} Soient $
(S{\stackrel{f}{\rightarrow}}T{\stackrel{g}{\rightarrow}}U)\in\overrightarrow{\mathbf{B}}^2=\overrightarrow{\mathbf{C}_{i_0}}^2$,
$\mu\in M$ et $a=(a_i)_{i\in I}\in S^\beta$ quelconques. Nous voulons vérifier que 
\[(g\circ f)^\beta_\mu(a)\subset (g^\beta_\mu \odot f^\beta_\mu) (a)
.\] Supposons $(g\circ f)^\beta_\mu(a)\neq\emptyset$, et soit $c=(c_i)_{i\in I}\in (g\circ f)^\beta_\mu(a)\subset U^\beta$, un élément quelconque de $(g\circ f)^\beta_\mu(a)$. 
Posons $t_0=\tau_{i_0}(a_{i_0})\in S^{\mathbf{h}_{i_0}}$, 
 $t_1=f^{\mathbf{h}_{i_0}}(t_0)\in T^{\mathbf{h}_{i_0}}$,
et $t_2=g^{\mathbf{h}_{i_0}}(t_1)\in U^{\mathbf{h}_{i_0}}$.
Par définition de $(g\circ f)^\beta_\mu(a)$, on a d'une part
\[
\tau_{i_0}(c_{i_0})
=(g\circ f)^{\mathbf{h}_{i_0}}(t_0)
=g^{\mathbf{h}_{i_0}}(f^{\mathbf{h}_{i_0}}(t_0))
=t_2,
\]
et d'autre part
\[\exists (\mathfrak{a}_i)_{i\in I}\in rb(R)^{-1}(\mu), \forall i \in I, \mathfrak{a}_i \rhd a_i,c_i.\]
Soit donc $(\mathfrak{a}_i)_{i\in I}$ une telle famille  de réalisations. Remarquons d'abord que ${\mathfrak{a}_{i_0} \rhd c_{i_0}} \Rightarrow {t_2=\tau_{i_0}(c_{i_0})\in D_{\mathfrak{a}_{i_0}}}$. Or, $t_1\leq_{\mathbf{h}_{i_0}}t_2$, d'où%
\footnote{Conformément à la remarque \ref{rmq domaine def des real} (page \pageref{rmq domaine def des real}).} 
 $t_1\in D_{\mathfrak{a}_{i_0}}$. De même, pour $i\neq {i_0}$, soit on a $\delta_i$  croissante, et dans ce cas on a $\delta_i(t_1)\leq_{\mathbf{h}_i} \delta_i(t_2)$, or $c_i=\mathfrak{a}_i(\tau_i(c_i))=\mathfrak{a}_i(\delta_i(t_2))$, d'où
 $\delta_i(t_1)\in D_{\mathfrak{a}_{i}}$,  soit on a $\delta_i$  décroissante, et dans ce cas on a $\delta_i(t_1)\leq_{\mathbf{h}_i} \delta_i(t_0)$, or $a_i=\mathfrak{a}_i(\tau_i(a_i))=\mathfrak{a}_i(\delta_i(t_0))$, d'où encore une fois $\delta_i(t_1)\in D_{\mathfrak{a}_{i}}$. L'hypothèse de monotonie de chaque $\delta_i$ permet ainsi de poser $b=(b_i)_{i\in I}$ avec $b_{i_0}=\mathfrak{a}_{i_0}(t_1)$ et, pour $i\neq {i_0}$, $b_i=\mathfrak{a}_i(\delta_i(t_1))$.
On a alors
\begin{itemize}
\item $t_1\in T^{\mathbf{h}_{i_0}}\Rightarrow b_{i_0}=\mathfrak{a}_{i_0} (t_1)\in   T^{\alpha_{i_0}}$,
et de même, pour $i\neq {i_0}$, $t_1\in T^{\mathbf{h}_{i_0}}\Rightarrow \delta_i(t_1)\in (\Delta_i T)^{\mathbf{h}_i}$, d'où $\mathfrak{a}_i(\delta_i(t_1))\in  (\Delta_i T)^{\alpha_i}$, autrement dit $b_i\in (\Delta_i T)^{\alpha_i}$, de sorte que 
\[b\in T^\beta,\]
\item par définition d'une réalisation,  $\tau_{i_0}(b_{i_0})=\tau_{i_0}(\mathfrak{a}_{i_0}(t_1))=t_1$ et, pour tout $i\neq {i_0}$,
\[
\tau_i(b_i)
=\tau_i(\mathfrak{a}_i(\delta_i(t_1)))
=\delta_i(t_1)
=\delta_i(\tau_{i_0}(b_{i_0})),
\]
\item enfin, par construction même,  $\mathfrak{a}_i\rhd a_i, b_i$ pour tout $i\in I$,
\end{itemize}
de sorte que, par définition de $f^\beta_\mu(a)$, on a $b\in f^\beta_\mu(a)$.\\

Par ailleurs, on a également
\begin{itemize}
\item $c\in U^\beta$,
\item $\tau_{i_0}(c_{i_0})=t_2=(g\circ f)^{\mathbf{h}_{i_0}}(t_0)=g^{\mathbf{h}_{i_0}}(t_1)=g^{\mathbf{h}_{i_0}}(\tau_{i_0}(b_{i_0}))$,
\item et, pour tout $i\in I$, $\mathfrak{a}_i\rhd b_i,c_i$,
\end{itemize}
de sorte que $c\in g^\beta_\mu(b)$.\\
Finalement, on a  $c\in \bigcup_{b\in f^\beta_\mu(a)}{g^\beta_\mu(b)}= (g^\beta_\mu \odot f^\beta_\mu) (a)$, d'où l'inclusion qu'il s'agissait de vérifier.
\begin{flushright}$\square$\end{flushright}

\begin{thm}[Stabilité sous-fonctorielle]\label{thm stabilite sous fonctorielle} Soit \[\mathcal{F}=(I,\mathcal{A}=(A_i)_{i\in I}, R, i_0, (\Delta_i,\delta_i)_{i\neq i_0})\] une famille interactive, de composantes les dynamiques sous-fonctorielles ouvertes efficientes
\[A_i=({\tau_i:(\alpha_i:\mathbf{C}_i\rightharpoondown \mathbf{P}^{\underrightarrow{\scriptstyle{L_i}}}})\looparrowright \mathbf{h}_i).\]
Les données suivantes :
\begin{itemize}
\item $\beta=\beta_\mathcal{F}$, 
la multi-dynamique sous-fonctorielle associée à $\mathcal{F}$ par la proposition \ref{prop multidynamique graphique est sous-fonctorielle},
\item $\mathbf{k}=\mathbf{h}_{i_0}$, l'horloge du chef d'orchestre $A_{i_0}$ de $\mathcal{F}$,
\item et $\rho:st(\beta)\rightarrow st(\mathbf{k})$, l'application définie par
\[\forall S\in\dot{\mathbf{C}}_{i_0}, \forall a=(a_i)_{i\in I}\in S^\beta, \rho(a)=\tau_{i_0}(a_{i_0}),\]
\end{itemize}
définissent une dynamique sous-fonctorielle ouverte
\[
[\mathcal{F}]_{\mathrm{p}}
=
(\rho:(\beta:\mathbf{C}_{i_0}\rightharpoondown \mathbf{P}^{\underrightarrow{\scriptstyle{M}}} )\looparrowright \mathbf{k}).
\]
\end{thm}
\paragraph{Preuve.}
Cela résulte immédiatement
\begin{itemize}
\item de la proposition \ref{prop multidynamique graphique est sous-fonctorielle},
\item du fait que pour tout $S\in \dot{\mathbf{C}}_{i_0}$ et tout $a\in S^\beta$, on a $\rho(a)\in S^\mathbf{k}$,
\item et du fait que
 que pour tout $(S{\stackrel{f}{\rightarrow}}T)\in\overrightarrow{\mathbf{C}_{i_0}}$, tout $\mu\in M$, tout $a=(a_i)_{i\in I}\in S^\beta$ et tout $b\in f^\beta_\mu(a)$,  on a
\[\rho(f^\beta_\mu(a))
=\tau_{i_0}(b_{i_0})
=f^{\mathbf{h}_{i_0}}(\tau_{{i_0}}(a_{i_0}))
=f^{\mathbf{k}}(\rho(a)). 
 \]
\end{itemize}
\begin{flushright}$\square$\end{flushright}

\subsubsection{Dynamique $[\mathcal{F}]_{\mathrm{p}}$ primo-engendrée par $\mathcal{F}$}

\begin{df}\label{df dynamique primo engendree}
La dynamique sous-fonctorielle ouverte $[\mathcal{F}]_{\mathrm{p}}$ associée dans le théorème \ref{thm stabilite sous fonctorielle}  à la famille interactive $\mathcal{F}$ est appelée la \emph{dynamique primo-engendrée par $\mathcal{F}$}.
\end{df}

\subsection{Diverses dynamiques engendrées par une famille interactive}\label{subs divers engendrements}

Comme indiqué dans la section \textbf{§\,4.2} de \cite{Dugowson:20150807}, 

\begin{quote}
[...] l'ensemble $M$ des paramètres de $[\mathcal{F}]_\mathrm{p}$ est en général \og trop gros\fg\, en ce sens que bien souvent le choix d'une valeur quelconque dans $M$ ne sera pas compatible avec le libre fonctionnement de la dynamique engendrée et, de ce fait, pourrait sembler peu naturel. 
\end{quote}

Afin de \emph{réduire} l'ensemble paramétrique, nous pouvons faire appel  à une relation d'équivalence $\sim$  sur $M$, choisie aussi judicieusement que possible,  pour former, conformément à la définition \ref{df reduc param multi} (page \pageref{df reduc param multi}), la dynamique sous-fonctorielle ouverte quotient $[\mathcal{F}]_\mathrm{p}/{\sim}$. En prenant pour relation $\sim$ la relation totale sur $M$, nous obtiendrons pour $[\mathcal{F}]_\mathrm{p}/{\sim}$ une mono-dynamique \og ouverte\fg, que nous noterons ${[\mathcal{F}]_\mathrm{m}}$, dont les transitions ne dépendent donc d'aucun paramètre. 
Dès lors il s'agit de
\begin{quote}
[...] trouver le juste équilibre entre l'ouverture excessive de ${[\mathcal{F}]_\mathrm{p}}$, offrant des paramètres en réalité souvent inutilisables, et la fermeture complète sur elle-même de ${[\mathcal{F}]_\mathrm{m}}$, qui n'offre plus aucune prise à l'interaction%
\footnote{Du moins si on se limites aux interactions normales, voir plus la haut la section \textbf{§\,\ref{subsubs interactions normales}}.
}  %
 avec d'autres dynamiques.\footnote{Extrait de la remarque 12 dans \cite{Dugowson:20150809}.}
\end{quote}

Dans la présente section \ref{subs divers engendrements}, reprenant la construction de la section \textbf{§\,4.2} de \cite{Dugowson:20150807}, nous considérons une classe particulière de relations d'équivalence sur l'ensemble paramétrique $M$, définies par le choix de \og tas paramétriques\fg\, au sein de chaque $L_i$, puis nous redonnons la définition, outre la dynamique primo-engendrée $[\mathcal{F}]_\mathrm{p}$ et  la dynamique mono-engendrée ${[\mathcal{F}]_\mathrm{m}}$ déjà citées,  de deux dynamiques engendrées de cette façon par la famille interactive $\mathcal{F}$.

\subsubsection{Tas paramétriques et équivalence sur $M$}\label{subsubsec tas param et equiv sur M}
Pour chaque $i\in I$, ayant fixé une certaine partie $N_i\subset L_i$  appelée le \emph{tas d'indice $i$}, partie intuitivement destinée à rassembler les valeurs du paramètre de la dynamique $A_i$ dont on considère qu'il appartient au libre fonctionnement de la dynamique engendrée que de les déterminer, et ayant ainsi constitué une famille $\mathcal{N}=(N_i)_{i\in I}$ de tas paramétriques, on considère sur $M\subset \Pi_I\mathcal{L}$ la relation d'équivalence $\sim_\mathcal{N}$ définie, pour tout couple $((\lambda_i)_{i\in I},(\lambda'_i)_{i\in I})\in M^2$, par
\[
(\lambda_i)_{i\in I}\sim_\mathcal{N}(\lambda'_i)_{i\in I}
\]
\[\Leftrightarrow\]
\[
\forall i\in I, \,
(\lambda_i=\lambda'_i)
\,\mathrm{ou}\,
(\lambda_i\in N_i\ni \lambda'_i).\]

La dynamique engendrée par $\mathcal{F}$ au sens des tas paramétriques $\mathcal{N}$ est alors\footnote{Voir la section \textbf{§\,\ref{subsubs quotient parametrique}}.}
\[[\mathcal{F}]_\mathcal{N}=[\mathcal{F}]_\mathrm{p}/\sim_\mathcal{N}.\]

\subsubsection{Dispositions $R$-compatibles} Pour préciser la manière dont les tas paramétriques sont définis dans les constructions des sections suivantes, nous aurons besoin de faire appel à la notion de \emph{disposition $R$-compatible} qui, dans le contexte mieux adapté
\footnote{Comme rappelé  en section \textbf{§\,\ref{subsubs relations binaires multiples}}, le \og transtypage\fg\, de $\mathcal{BM}_I$ à  $\mathcal{R}_{2I}$ est réalisé par l'application $rd$ qui, en particulier, transforme canoniquement toute relation binaire multiple $R\in\mathcal{BM}_{(\mathcal{Z}^*,\mathcal{L})}$ indexée sur $I$ en une relation multiple $rd(R)$ d'index $2I=I\sqcup I$.}
 des relations multiples
indexées par $2I=I\sqcup I=I\times\{0,1\}$ plutôt que dans celui de relations binaires multiples indexées par $I$, généralise aux configurations partielles la notion de $R$-configuration.

\begin{df}[Disposition] 
On appelle \emph{disposition entière pour $\mathcal{A}$} tout élément $q$, indexé par $2I$,  du produit\footnote{Sur l'ordre des facteurs, voir la remarque \ref{rmq ordre des facteurs}.} $(\Pi {\mathcal{Z}^*})\times(\Pi {\mathcal{L}})\simeq \Pi \mathcal{E}$, qui soit de la forme  $q=rd_\mathcal{E}(\varsigma)$, avec $\varsigma\in \vert \Omega_\mathcal{A} \vert\subset \Pi \mathcal{E}$. 
Plus généralement, pour toute partie  $W\subset 2I$, on appelle $W$-\emph{disposition pour $\mathcal{A}$} toute famille $q=(q_w)_{w\in W}$ d'éléments pris  --- selon que $w$ est, dans $2I$, respectivement de la forme $(i,0)$ ou $(i,1)$ ---  dans les ensembles constituant respectivement la famille $\mathcal{Z}^*$ ou ceux de la famille $\mathcal{L}$, telle qu'il existe une disposition entière $rd_\mathcal{E}(\varsigma)$ dont la restriction à $W$ soit égale à $(q_w)_{w\in W}$ :
\[\exists \varsigma\in\vert\Omega_\mathcal{A}\vert,\forall w\in W, q_w=rd_\mathcal{E}(\varsigma)_w.\] On appelle enfin \emph{disposition pour $\mathcal{A}$} toute  $W$-disposition pour un certain $W\subset 2I$.
\end{df}

Intuitivement, une disposition  n'est rien d'autre, à ordre des facteurs près, qu'une configuration partielle, autrement dit une famille constituée de réalisations non vides et de valeurs paramétriques cohérentes pour $\mathcal{A}$.\\

Pour tout $W\subset 2I$, nous noterons $\Pi_W{[\mathcal{A}]}$ l'ensemble des $W$-dispositions pour $\mathcal{A}$. L'ensemble des dispositions pour $\mathcal{A}$ s'écrit ainsi $\bigcup_{W\subset 2I}{\Pi_W{[\mathcal{A}]}}$, 
et on a également $\Pi_{2I}{[\mathcal{A}]}\simeq\Pi \mathcal{E}$, $\Pi_{I\times\{0\}}{[\mathcal{A}]}\simeq\Pi  \mathcal{Z}^*$ et $\Pi_{I\times\{1\}}{[\mathcal{A}]}\simeq\Pi  \mathcal{L}$.

\begin{df}[Dispositions $R$-compatibles]
Soit $W\subset 2I$. Une $W$-disposition  $q=(q_w)_{w\in W}$  est dite \emph{$R$-compatible} (ou \emph{compatible} avec $R$) si elle est la restriction à $W$ d'une disposition entière $rd_\mathcal{E}(\varsigma)\in \vert rd(R)\vert=rd(\vert R\vert)$. Autrement dit,  $q$  est  $R$-compatible si
\[\exists \varsigma\in\vert R\vert,\forall w\in W, q_w=rd_\mathcal{E}(\varsigma)_w.\] Nous écrirons \[q\wr R\] pour exprimer que $q$ est une disposition $R$-compatible.
\end{df}

Par exemple, on a $(L_k\ni\lambda_k \wr R) \Leftrightarrow (\exists\mu\in M=Im(rb(R)),\lambda_k=\mu_k)$.\\

Par ailleurs, étant donnés $X\subset 2I$, $Y\subset 2I$  et deux familles $q=(q_x)_{x\in X}\in \Pi_X\mathcal{E}$ et $r=(r_y)_{y\in Y}\in \Pi_Y\mathcal{E}$  \emph{compatibles entre elles} au sens où pour tout $w\in X\cap Y$ on a $q_w=r_w$ --- ce qui est le cas notamment si $X\cap Y=\emptyset$ --- nous noterons $q\vee r$ la famille $q\vee r=s=(s_w)_{w\in X\cup Y}$ telle que, pour tout $w\in X\cup Y$, on a $w\in X \Rightarrow s_w=q_w$ et $w\in Y \Rightarrow s_w=r_w$. Plus généralement, on définit l'opération $\vee$ pour les familles de dispositions deux à deux compatibles et, sous réserve de compatibilité, $\vee$ est associative, commutative, et admet la disposition vide pour élément neutre.  
Nous utiliserons en particulier cette notation avec des familles de la forme \[(\mathfrak{a}_j)_{j\in J\subset I}\in \Pi_J(\mathcal{Z}^*),\] ou de la forme \[(\lambda_k)_{k\in K\subset I}\in \Pi_K(\mathcal{L}).\] 
Par exemple, écrire que l'on a $\lambda=(\lambda_i)_{i\in I}\in M$, ce qui revient à écrire $\lambda\wr R$, équivaut encore à
\[\exists \mathfrak{a}\in \Pi_I(\mathcal{Z}^*), \lambda \vee \mathfrak{a}\in rd(\vert R\vert). \]

\subsubsection{Dynamique $[\mathcal{F}]_\mathrm{f}$ fonctionnellement engendrée par $\mathcal{F}$}\label{subsubsec dyna fonct eng}

Pour tout $k\in I$, on définit le tas fonctionnel $N^{\mathrm{f}}_k\subset L_k$ de la façon suivante : un élément $l_k\in L_k$ vérifie $l_k\in N^{\mathrm{f}}_k$ si et seulement si $l_k$ est $R$-compatible\footnote{\label{footnote condition l compatible}Logiquement, la condition \og $l_k$ est $R$-compatible\fg\, ne change rien, puisqu'au bout du compte la relation d'équivalence définie par les tas le sera sur l'ensemble $M=Im(rb(R))$, mais c'est sans doute plus clair de se limiter aux $l_k$ effectivement concernés.} et si on a
\[
\forall \mathfrak{a}\in\Pi_{i\neq k}{\mathcal{S}^*_{A_i}},
\forall \lambda_k\in L_k,
\left(
\left.
\begin{tabular}{c}
$(\mathfrak{a}\vee  l_k) \wr R$\\ 
$(\mathfrak{a}\vee  \lambda_k)\wr R$\\ 
\end{tabular} 
\right\rbrace
\Rightarrow
l_k=\lambda_k
\right).
\]

Intuitivement, on place dans le tas $N^{\mathrm{f}}_k$ les paramètres de la dynamique ouverte $A_k$ dont la valeur est déterminée \emph{via} l'interaction $R$ par les réalisations des \emph{autres} dynamiques en jeu.

La famille $\mathcal{N}^{\mathrm{f}}=(N^{\mathrm{f}}_k)_{k\in I}$ des tas paramétriques fonctionnels étant ainsi définie, on applique le procédé décrit dans la section \ref{subsubsec tas param et equiv sur M}, obtenant ainsi une relation d'équivalence $\sim_{\mathrm{f}}$ sur $M$, et on appelle \emph{dynamique ouverte sous-fonctorielle fonctionnellement engendrée par la famille interactive $\mathcal{F}$}, et l'on note $[\mathcal{F}]_\mathrm{f}$, la dynamique sous-fonctorielle ouverte
\[ [\mathcal{F}]_\mathrm{f}=[\mathcal{F}]_\mathrm{p}/\sim_{\mathrm{f}},\] d'ensemble de paramètres $M/\sim_{\mathrm{f}}$.

\subsubsection{Dynamique $[\mathcal{F}]_\mathrm{s}$ souplement engendrée par $\mathcal{F}$}\label{subsubsec dyna soupl eng}

\begin{rmq}
La définition donnée ci-après de la dynamique que nous diront \emph{souplement engendrée} par une famille interactive nous a été suggéré par l'examen de ce qui devrait être considéré comme libre paramètre ou non dans le cas d'un ressort soumis à différents jeux de contraintes tels que la connaissance du comportement du ressort permette de savoir à quel jeu de contrainte il est soumis. Nous laissons ici  au lecteur une telle étude en exercice. 
\end{rmq}

Pour tout $k\in I$, on définit le tas fonctionnel $N^{\mathrm{s}}_k\subset L_k$ comme l'ensemble des éléments bloqués de $L_k$, un élément de $L_k$ étant dit \emph{bloqué} (pour l'interaction $R$) s'il n'est pas libre, tandis qu'un élément $\lambda_k\in L_k$ est dit \emph{libre} ou \emph{souple} s'il est $R$-compatible\footnote{La remarque de la note \ref{footnote condition l compatible} ci-dessus s'applique encore ici, et s'appliquerait aussi bien aux éléments bloqués.} et si quel que soit $\mathfrak{a}_k\in\mathcal{S}^*_{A_k}$, quel que soit $\mu\in \Pi_{j\neq k}L_j$ et quel que soit $\mathfrak{b}\in \Pi_{j\neq k}\mathcal{S}^*_{A_j}$ on a l'implication
\[
\left.
\begin{tabular}{c}
$(\lambda_k\vee \mathfrak{a}_k\vee \mu)\wr R$\\ 
$(\mu\vee \mathfrak{b})\wr R$\\ 
\end{tabular} 
\right\rbrace
\Rightarrow
(\lambda_k\vee \mu\vee \mathfrak{a}_k\vee \mathfrak{b}) \wr R.
\]

La famille $\mathcal{N}^{\mathrm{s}}=(N^{\mathrm{s}}_k)_{k\in I}$ constituée des tas de paramètres bloqués étant ainsi définie, on applique le procédé décrit dans la section \ref{subsubsec tas param et equiv sur M}, obtenant ainsi une relation d'équivalence $\sim_{\mathrm{s}}$ sur $M$, et on appelle \emph{dynamique ouverte sous-fonctorielle souplement engendrée par la famille interactive $\mathcal{F}$}, et l'on note $[\mathcal{F}]_\mathrm{s}$, la dynamique sous-fonctorielle ouverte
\[ [\mathcal{F}]_\mathrm{s}=[\mathcal{F}]_\mathrm{p}/\sim_{\mathrm{s}},\] d'ensemble de paramètres $M/\sim_{\mathrm{s}}$.

Intuitivement, sont mis dans les tas de paramètres bloqués ceux dont le choix par un agent extérieur à la famille interactive considérée pourrait être en retour remis en cause par ce que nous pourrions appeler le libre fonctionnement de cette famille, représenté ici par $\mathfrak{a}_k$, $\mu\in \Pi_{j\neq k}L_j$ et  $\mathfrak{b}$. 

\subsubsection{Dynamique $[\mathcal{F}]_\mathrm{m}$ mono-engendrée par $\mathcal{F}$}\label{subsubsec dyna mono eng}

Prenant pour relation d'équivalence sur $M$ la relation d'équivalence maximale $\sim_{\mathrm{m}}$, de sorte que $M/\sim_{\mathrm{m}}$ est réduit à un point, on obtient comme annoncé précédemment la mono-dynamique sous-fonctorielle  engendrée par la famille interactive $\mathcal{F}$,  encore appelée \emph{dynamique sous-fonctorielle  mono-engendrée par $\mathcal{F}$} :
\[ [\mathcal{F}]_\mathrm{m}=[\mathcal{F}]_\mathrm{p}/\sim_{\mathrm{m}}.\] 
 Bien entendu, cette mono-dynamique scandée peut toujours être vue comme une dynamique ouverte, l'ensemble des valeurs prises par le paramètre se réduisant à un singleton. Cette construction revient à prendre pour tas paramétriques d'indice $k$ l'ensemble $L_k$ lui-même ou, de façon équivalente\footnote{Comme indiqué dans la note \ref{footnote condition l compatible}.},  l'ensemble des valeurs de $L_k$ qui sont $R$-compatibles.
 
\subsection{Exemples d'engendrement dynamique}\label{subs exemple engendrement}

Souhaitant conserver à cet article une taille raisonnable, nous ne donnerons ici que deux exemples, même s'il est clair que de nombreux autres exemples devraient être développés pour éclairer ne serait-ce que les principes fondamentaux de la théorie du dynamisme sous-fonctoriel.

\subsubsection{Exemple de la famille interactive canoniquement associée à une dynamique ouverte}

\begin{exm}
 \`{A} toute dynamique sous-fonctorielle ouverte efficiente 
\[A=({\tau:(\alpha:\mathbf{C}\rightharpoondown \mathbf{P}^{\underrightarrow{\scriptstyle{L}}}})\looparrowright \mathbf{h})\]
nous associons la famille interactive $\mathcal{F}=[A]$ définie par
\[[A]=(I=\{0\},\mathcal{A}=(A), R=\Omega_\mathcal{A},i_0=0, (\Delta_i,\delta_i)=\emptyset).\]
On vérifie sans peine que la dynamique primo-engendrée $[\mathcal{F}]_\mathrm{p}=[[A]]_\mathrm{p}=B$ est de la forme
\[B=({\tau:(\beta:\mathbf{C}\rightharpoondown \mathbf{P}^{\underrightarrow{\scriptstyle{M}}}})\looparrowright \mathbf{h})\]
avec 
\begin{itemize}
\item $M=\{\mu\in L,\mathcal{S}^*_{(A,\mu)}\neq\emptyset\}\subset L$,
\item pour tout $S\in\dot{\mathbf{C}}$, $S^\beta=S^\alpha$,
\item pour toute arête $(e:S\rightarrow T)\in \overrightarrow{\mathbf{C}}$, tout état  $a\in S^\alpha$ et tout paramètre $\mu\in M$,
\[ 
e^\beta_\mu(a)=
\{
b\in T^\alpha,  
(\tau (b)=e^{\mathbf{h}}(\tau(a)))
\,\,\mathrm{et}\,\,
(\exists \mathfrak{a}\in \mathcal{S}^*_{(A,\mu)}, 
\mathfrak{a}\triangleright a, b)
\}. 
\]
\end{itemize}

De plus, on a $[[A]]_\mathrm{p}=[[A]]_\mathrm{f}=[[A]]_\mathrm{s}$, tandis que la dynamique mono-engendrée $C=[[A]]_\mathrm{m}$ s'écrit 
\[C=({\tau:(\gamma:\mathbf{C}\rightharpoondown \mathbf{P}})\looparrowright \mathbf{h})\]
avec 
\begin{itemize}
\item pour tout $S\in\dot{\mathbf{C}}$, $S^\gamma=S^\alpha$,
\item pour toute arête $(e:S\rightarrow T)\in \overrightarrow{\mathbf{C}}$ et tout état  $a\in S^\alpha$ 
\[ 
e^\gamma(a)=
\{
b\in T^\alpha,  
(\tau (b)=e^{\mathbf{h}}(\tau(a)))
\,\,\mathrm{et}\,\,
(\exists \mathfrak{a}\in \mathcal{S}^*_{A}, 
\mathfrak{a}\triangleright a, b)
\}. 
\]
\end{itemize}

\begin{rmq}
Tandis que les réalisations de $A$ dépendent des possibilités données par ses transitions, les transitions de la dynamique $[[A]]_\mathrm{p}$  dépendent des réa\-li\-sations de $A$. L'application $A\mapsto [[A]]_\mathrm{p}$ met ainsi en évidence la dialectique entre  les \og possibilités de principe\fg\, qu'indiquent  les transitions et les \og possibilités effectives\fg qu'expriment les réalisations, dialectique qui a fait l'objet de notre exposé filmé \cite{Dugowson:20150506}. Bien entendu, il y a une large classe de dynamiques, que nous dirons  \emph{régulières}, pour lesquelles $A=[[A]]_\mathrm{p}$. Nous donnons ci-après deux exemples pour lesquels ce n'est pas le cas.
\end{rmq}

\begin{exm}[Une dynamique déterministe sur $\mathbf{R}_+$ non régulière]\label{exm une dynamique sur R non reguliere}
Prenons \[A=(\tau:(\alpha:\mathbf{C}\rightarrow \mathbf{Sets})\looparrowright \mathbf{h})\] 
avec 
\begin{itemize}
\item $\mathbf{C}=(\mathbf{R}_+,+)$, le monoïde des réels positifs dont l'unique objet sera noté $\bullet$,
\item pour ensemble d'états $st(\alpha)=\bullet^\alpha=\mathbf{R}\sqcup ({\mathbf{R}_+}\times \{1\})$,
\item pour horloge : $\bullet^\mathbf{h}=\mathbf{R}$, avec $\forall d\in\mathbf{R}_+, \forall r\in\mathbf{R}, d^\mathbf{h}(r)=r+d$,
\item pour scansion : $\forall r\in \mathbf{R}, \tau(r)=r$ et $\forall r\in \mathbf{R}_+, \tau(r,1)=r$,
\item pour transitions $d^\alpha$ associée à $d\in\mathbf{R}_+$ : $\forall r\in \mathbf{R}, d^\alpha(r)=r+d$ et $\forall r\in \mathbf{R}_+, \tau(r,1)=(r+d,1)$.
\end{itemize}

$A$ est une dynamique déterministe, donc fonctorielle. Cependant, l'horloge choisie interdit l'existence de réalisations passant par un état quelconque pris dans ${\mathbf{R}_+}\times \{1\}$, de sorte que pour $[[A]]_\mathrm{p}$ ces états sont hors-jeu, et la dynamique $[[A]]_\mathrm{p}$, quasi-déterministe mais non déterministe, n'est donc pas fonctorielle, et $[[A]]_\mathrm{p}\neq A$. Par contre, $[[A]]_\mathrm{p}=(\tau:\beta\looparrowright \mathbf{h})$ vérifie encore $e^\beta\odot d^\beta=(d+e)^\beta$. L'exemple suivant montre que cette relation elle-même n'est plus nécessairement satisfaite, même avec une dynamique $A$ déterministe.
\end{exm} 

\begin{exm}[Autre exemple de dynamique déterministe non régulière]\label{exm autre dynamique non reguliere}
Soit \[A=(\tau:(\alpha:\mathbf{C}\rightarrow \mathbf{Sets})\looparrowright \mathbf{h})\] l'unique mono-dynamique ouverte déterministe telle que
\begin{itemize}
\item $\mathbf{C}$ est la catégorie à quatre objets $\dot{\mathbf{C}}=\{S,U,V,T\}$ et dix flèches $\overrightarrow{\mathbf{C}}=\{Id_S, \overrightarrow{SU}, \overrightarrow{SV}, \overrightarrow{SUT}, \overrightarrow{SVT}, Id_U, \overrightarrow{UT}, Id_V, \overrightarrow{VT}, Id_T\}$  où, mises à part les identités, la première lettre du nom d'une flèche désigne son domaine et la dernière lettre son codomaine (par exemple, les flèches de domaine $S$ sont $Id_S$, $\overrightarrow{SU}$, $\overrightarrow{SV}$, $\overrightarrow{SUT}$ et $\overrightarrow{SVT}$, celles de codomaine $T$ sont $\overrightarrow{SUT}$, $\overrightarrow{SVT}$, $\overrightarrow{UT}$,  $\overrightarrow{VT}$, $Id_T$, et ainsi de suite) et où la composition se fait selon ce que ces noms suggère, par exemple $\overrightarrow{VT}\circ \overrightarrow{SV} = \overrightarrow{SVT}\neq \overrightarrow{SUT}=\overrightarrow{UT}\circ \overrightarrow{SU}$, etc.
\item les ensembles d'états sont $S^\alpha=\{s, s'\}$, $U^\alpha=\{u, u'\}$, $V^\alpha=\{v\}$, $T^\alpha=\{t, t'\}$, ces sept états étant deux à deux distincts,
\item l'horloge $\mathbf{h}$ est l'horloge essentielle%
\footnote{Voir l'exemple \ref{exm horloge essentielle et existentielle} page \pageref{exm horloge essentielle et existentielle}.} 
de $\mathbf{C}$, $\mathbf{h}=\zeta_\mathbf{C}$, pour laquelle $S^\mathbf{h}=\{S\}$, $U^\mathbf{h}=\{U\}$, etc., et qui vérifie notamment $\overrightarrow{SUT}^\mathbf{h}=\overrightarrow{SVT}^\mathbf{h}$ puisque ces deux expressions désignent toutes deux l'unique application $\{S\}\rightarrow \{T\}$,
\item admettant notamment les transitions suivantes, les autres découlant de la nature déterministe (et en particulier fonctorielle) de cette dynamique : $\overrightarrow{SU}^\alpha(s)=u$, $\overrightarrow{SU}^\alpha(s')=u'$, $\overrightarrow{UT}^\alpha(u)=t$, $\overrightarrow{UT}^\alpha(u')=t'$, $\overrightarrow{VT}^\alpha(v)=t$,
\end{itemize}
et soit $(\tau:(\beta \looparrowright \mathbf{h}))=[[A]]_\mathrm{p}=[[A]]_\mathrm{m}$ la dynamique engendrée. 
Alors on vérifie facilement que, s'il y a des réalisations de $A$ qui passent par $s'$ et $v$, et qu'il y en a une qui passe par $v$ et $t$, aucune ne passe par $s'$ et $t$, ni d'ailleurs par $s'$ et $t'$, de sorte que l'on a $\overrightarrow{SVT}^\beta(s')=\emptyset\subsetneqq \{t\}=(\overrightarrow{VT}^\beta\odot \overrightarrow{SV}^\beta)(s')$.
\end{exm} 
\end{exm}

\subsubsection{Exemple de la dynamique $\mathbb{S}= $ \textcjheb{/s} $ =[$\textcjheb{why}$]_\mathrm{m}=[\mathbb{WHY}]_\mathrm{m}$}\label{subsubs exemple shin}

En utilisant les données de la section \textbf{§\,\ref{subsubs famille why}}, on peut montrer d'abord\footnote{Pour cet exemple, toutes les vérifications  sont laissées en exercice au lecteur.} que la dynamique {$\mathbb{T}=[\textcjheb{why}]_\mathrm{p}=[\mathbb{WHY}]_\mathrm{p}=(\tau:\beta\looparrowright \mathbf{h}_\mathbb{T})$} primo-engendrée par la famille interactive \textcjheb{why}$=\mathbb{WHY}$ admet
\begin{itemize}
\item pour moteur $\mathbf{C}_{\mathbb{T}}=\mathbf{C}_{\textcjheb{h}}=(\mathbf{R}_+,+)$,
\item pour horloge $\mathbf{h}_\mathbb{T}=\mathbf{h}_{\textcjheb{h}}$, 
de sorte que pour tout $t\in st(\mathbf{h}_\mathbb{T})=]0,+\infty[$ 
et  tout $r\in\mathbf{R}_+$, on a
$r^{\mathbf{h}_\mathbb{T}}(t)=t+r$,
\item pour ensemble paramétrique $M$ l'ensemble des $(\omega,\gamma,*)\in \mathcal{C}^1_{\vartriangleright}\times Lip^1\times\{*\}$ tels qu'il existe un triplet $(t_\omega,\overline{s}_\gamma, t_\gamma)\in \overline{\mathbf{R}}_+^3$ avec $0<t_\omega\leq \overline{s}_\gamma\leq t_\gamma\leq +\infty$ et une fonction $\lambda\in\mathcal{C}(]-\infty,t_\gamma[)$ qui vérifient
\subitem - $\omega=\lambda_{\vert ]-\infty, t_\omega[}\in \mathcal{C}^1_{\vartriangleright}$,
\subitem - $\gamma=\lambda_{\vert ]0, t_\gamma[}\in Lip^1$,
\subitem - $\overline{s}_\gamma=\sup \{s\in]0,t_\gamma], \lambda_{\vert ]-\infty, s[}\in \mathcal{C}^1\}$,
\item pour ensemble d'états $st([\mathbb{WHY}]_\mathrm{p})$ l'ensemble des quadruplets $(t,r,f,w)$ avec $t\in\mathbf{R}_+^*$, $r\in \mathbf{R}$, $f\in\mathcal{C}^1(]-\infty,t[)$ et $w\in\mathcal{C}$, 
\item pour scansion  $(t,r,f,w)\mapsto\tau(t,r,f,w)=t$,
\item pour transition de paramètre $(\omega,\gamma,*)\in M$ associée à $d\in\mathbf{R}_+$ la fonction $(t,r,f,w)\mapsto d^\beta_{(\omega,\gamma,*)}(t,r,f,w)=(t',r',f',w')$ qui --- si $0<t\leq t+d\leq t_\omega$,  $\lim_{s\rightarrow t^-}f(s)=r$,
 $f=\omega_{\vert ]-\infty,t[}$ et $\omega \rightY w$ (qui sont les conditions pour rester dans le jeu) --- est donnée par
\[
d^\beta_{(\omega,\gamma,*)}(t,r,f,w)
=
(t'=t+d, r'=\omega(t+d)=\gamma(t+d), f'=\omega_{]-\infty,t+d[},w'=w),
\]
\end{itemize}
où, le cas échéant, on prolonge par continuité $\omega$ en $t_\omega$ et $\gamma$ en $t_\gamma$.\\

On vérifie ensuite que $[\mathbb{WHY}]_\mathrm{p}=[\mathbb{WHY}]_\mathrm{f}=[\mathbb{WHY}]_\mathrm{s}$. \\

Posant $\mathbb{S}= \textcjheb{/s} = [\textcjheb{why}]_\mathrm{m}=[\mathbb{WHY}]_\mathrm{m}=(\tau:\varphi\looparrowright\mathbf{h}_\mathbb{S})$ --- dynamique dont par définition le moteur, l'horloge, les états et la scansion sont les mêmes que ceux de $\mathbb{T}$ --- on vérifie enfin que pour tout état $(t,r,f,w)$ qui est dans le jeu --- c'est-à-dire tel que $\lim_{s\rightarrow t^-}f(s)=r$, $f_{\vert [0,t[}\in Lip^1([0,t[)$ et $f\rightY w$ --- et pour tout $d\geq 0$, l'ensemble $d^\varphi(t,r,f,w)$ est constitué des $(t',r',f',w')\in st(\mathbb{S})$ qui vérifient

\begin{itemize}
\item $t'=t+d$,
\item $r'\in[r-d,r+d]$,
\item $(f')_{\vert ]-\infty,t[}= f_{\vert ]-\infty,t[}$,
\item $(f')_{\vert ]0,t+d[}\in Lip^1$,
\item $\lim_{s\rightarrow {t'}^-}f'(s)=r'$
\item $f'\rightY w'=w$.
\end{itemize}

\bibliographystyle{plain}




\tableofcontents

\end{document}